\documentclass[12pt]{amsart}
\usepackage{latexsym, amsmath, amssymb, amsthm}
\usepackage{epsfig}
\usepackage{amscd}
\usepackage{latexsym}
\usepackage{pstricks,pst-node,cite}

\makeatletter
\newtheorem*{rep@theorem}{\rep@title}
\newcommand{\newreptheorem}[2]{%
\newenvironment{rep#1}[1]{%
 \def\rep@title{#2 \ref{##1}}%
 \begin{rep@theorem}}%
 {\end{rep@theorem}}}
\makeatother

\newtheorem{theorem}{Theorem}[section]
\newreptheorem{theorem}{Theorem}

\newtheorem{cor}[theorem]{Corollary}

\newtheorem{conj}[theorem]{Conjecture}
\theoremstyle{remark}

\newcommand{\tndi}{tndi_{\sum}}

\newcounter{fignum}
\ifx\blackandwhite\undefined
  \newrgbcolor{darkred}{.9 0 0}\newrgbcolor{emgreen}{0 .9 0}
  \newrgbcolor{pup}{.7  0 .9}\newrgbcolor{xxx}{.8 .8 .8}
  \newrgbcolor{lightgray}{.99 .99 0}\newrgbcolor{darkgray}{.8 .8 .3} \else
\fi
\psset{linecolor=blue}

\psset{linewidth=.7pt,dash=3pt 3pt,doublesep=.05,dotsize=1pt 5}
\SpecialCoor

\begin{document}

\title[A proper total coloring distinguishing adjacent vertices by sums]{A proper total coloring distinguishing adjacent vertices by sums of some product graphs}

\author{Hana Choi}
\address{Department of Mathematics \\Kyonggi University
\\ Suwon, 443-760 Korea}

\author{Dongseok Kim}
\address{Department of Mathematics \\Kyonggi University
\\ Suwon, 443-760 Korea}
\email{dongseok@kgu.ac.kr}

\author{Sungjin Lee}
\address{Department of Mathematics \\Kyonggi University
\\ Suwon, 443-760 Korea}

\author{Yeonhee Lee}
\address{Department of Mathematics \\Kyonggi University
\\ Suwon, 443-760 Korea}

\begin{abstract}
In this article, we consider a proper total coloring distinguishes adjacent vertices by sums,
if every two adjacent vertices have different total sum of colors of the
edges incident to the vertex and the color of the vertex. Pilsniak and Wozniak~\cite{PW}
first introduced this coloring and made a conjecture that the minimal
number of colors need to have a proper total coloring distinguishes adjacent vertices by sums is less
than or equal to the maximum degree plus $3$.
We study proper total colorings distinguishing adjacent vertices by sums of some graphs and their products.
We prove that these graphs satisfy the conjecture.
\end{abstract}


\maketitle

\section{Introduction}

Let $\Gamma$ be a finite simple graph with vertex set $V(\Gamma)$ and
edge set $E(\Gamma)$. \emph{An $n$ coloring} $\phi$ of $\Gamma$ is a function from
$V(\Gamma)$ to $\{ 1, 2, \ldots, n \}$. A coloring $\phi$ is \emph{proper} if
$\phi(u) \not= \phi(v)$ for all edges $\{u,v\} \in E(\Gamma)$. \emph{The chromatic number}
$\chi(\Gamma)$ of a graph $\Gamma$ is the smallest number of colors needed
to color properly the vertices of $\Gamma$. Since the exploratory paper by Dirac \cite{dirac},
the chromatic number has been in the center of graph theory
research. Its rich history can be found in several articles
\cite{HM, thmas}.

This concept naturally expand to \emph{an edge coloring} which is a function from
$E(\Gamma)$ to $\{ 1, 2, \ldots, n \}$. An edge coloring $\phi$ is \emph{proper} if
$\phi(e_1) \not= \phi(e_2)$ if $e_1$ and $e_2$ have a common vertex.
The minimum required number of colors for a proper edge coloring of
a given graph is called \emph{the chromatic index} of the graph, denoted by $\chi'(\Gamma)$.
A very famous theorem by Vizing~\cite{Vizing} slanted that
for simple graph, the chromatic index of the graph is either its maximum degree $\Delta$ or $\Delta+1$.
For some graphs, such as bipartite graphs and high-degree planar graphs,
the chromatic index is always $\Delta$.

\emph{A total coloring} is a coloring on the vertices and edges of
a graph $\Gamma$ such that no adjacent vertices have the same color,
no adjacent edges have the same color and no edge and its end-vertices
are assigned the same color.
The \emph{total chromatic number} $\chi''(\Gamma)$ of a graph $\Gamma$
is the least number of colors needed in a total coloring of $\Gamma$.
Some properties of total chromatic number $\chi''(\Gamma)$ are as follows;
$\chi''(\Gamma)\ge \Delta +1$ and for upper bound for $\chi''(\Gamma)$,
Molloy and Reed~\cite{MR} first found that $\chi''(\Gamma)\le \Delta +10^{26}$
and $\chi''(\Gamma)\le ch'(\Gamma) + 2$, where $ch'(\Gamma)$
is \emph{the edge choosability}
which is the least number $k$ such that every instance of the list
edge-coloring that has $\Gamma$ as
its underlying graph and that provides at least $k$ allowed
colors for each edge of $\Gamma$ has a proper coloring.
A long standing open problem about the total coloring is
$\chi''(\Gamma)\le \Delta +2$ arose by Behzad and Vizing~\cite{JT}.
We refer to \cite{JT} for more problems on vertex, edge and total colorings.

For total coloring $\phi$ on a simple graph $\Gamma$, \emph{the color set of a vertex} $v$ is
defined by $C(v)=\{ \phi(u), \phi(\{u,v\}) | \{u,v\} \in E(\Gamma)\}$.
Z. Zhang et al.~\cite{ZCLYLW} introduced a new concept that a total coloring of a
graph $\Gamma$ is \emph{an adjacent vertex distinguishing
total coloring (AVD total coloring)}
if $C(u) \neq C(v)$ for all $\{u,v\} \in E(\Gamma)$.
\emph{The adjacent-vertex-distinguishing-total-chromatic number},
denoted by $\chi_{at}(\Gamma)$ of a graph $\Gamma$ is the
least number of colors needed in an AVD-total-coloring of $\Gamma$.
There have been several articles studying AVD-total-coloring of graphs~\cite{CG, CZ, Hulgan, HWY, MP, WH, WW}

\begin{figure}
$$
\begin{pspicture}[shift=-2.5](-3,-3)(3,2.7)
\psline(2;90)(2;18)(2;-54)(2;-126)(2;-198)(2;90)
\psline(2;90)(2;-54)(2;-198)(2;18)(2;-126)(2;90)
\pscircle[fillcolor=lightgray, fillstyle=solid, linewidth=1pt](2;90){.3}
\pscircle[fillcolor=lightgray, fillstyle=solid, linewidth=1pt](2;18){.3}
\pscircle[fillcolor=lightgray, fillstyle=solid, linewidth=1pt](2;-54){.3}
\pscircle[fillcolor=lightgray, fillstyle=solid, linewidth=1pt](2;-126){.3}
\pscircle[fillcolor=lightgray, fillstyle=solid, linewidth=1pt](2;-198){.3}
\pscircle[fillcolor=white, linecolor=white, fillstyle=solid, linewidth=1pt](.65;90){.15}
\pscircle[fillcolor=white, linecolor=white, fillstyle=solid, linewidth=1pt](.65;18){.18}
\pscircle[fillcolor=white, linecolor=white, fillstyle=solid, linewidth=1pt](.65;-54){.15}
\pscircle[fillcolor=white, linecolor=white, fillstyle=solid, linewidth=1pt](.65;-126){.18}
\pscircle[fillcolor=white, linecolor=white, fillstyle=solid, linewidth=1pt](.65;-198){.18}
\pscircle[fillcolor=white, linecolor=white, fillstyle=solid, linewidth=1pt](1.65;126){.15}
\pscircle[fillcolor=white, linecolor=white, fillstyle=solid, linewidth=1pt](1.65;54){.18}
\pscircle[fillcolor=white, linecolor=white, fillstyle=solid, linewidth=1pt](1.65;-18){.18}
\pscircle[fillcolor=white, linecolor=white, fillstyle=solid, linewidth=1pt](1.65;-90){.15}
\pscircle[fillcolor=white, linecolor=white, fillstyle=solid, linewidth=1pt](1.65;-166){.18}
\rput(2;90){{$3$}}\rput(2;18){{$5$}}\rput(2;-54){{$7$}}
\rput(2;-126){{$2$}}\rput(2;-198){{$4$}}
\rput(.65;90){{\footnotesize{1}}}\rput(.65;18){{\footnotesize{5}}}
\rput(.65;-54){{\footnotesize{7}}}
\rput(.65;-126){{\footnotesize{2}}}\rput(.65;-198){{\footnotesize{6}}}
\rput(1.65;126){{\footnotesize{7}}}\rput(1.65;54){{\footnotesize{4}}}
\rput(1.65;-18){{\footnotesize{6}}}
\rput(1.65;-90){{\footnotesize{1}}}\rput(1.65;-166){{\footnotesize{3}}}
\rput(2.6;90){{$\{1,2\}^{c}$}}\rput(2.9;20){{$\{2,3\}^{c}$}}
\rput(2.9;-52){{$\{3,4\}^{c}$}}\rput(2.9;-128){{$\{4,5\}^{c}$}}
\rput(2.9;-200){{$\{5,6\}^{c}$}} \rput(0,-2.9){{$(a)$}}
\end{pspicture}
\quad\quad
\begin{pspicture}[shift=-2.5](-3,-3)(3,2.7)
\psline(0,0)(2;0)\psline(0,0)(2;60)\psline(0,0)(2;120)\psline(0,0)(2;180)
\psline(0,0)(2;240)\psline(0,0)(2;300)
\psline(1,1.72)(-1,1.72)(-2,0)(-1,-1.72)(1,-1.72)(2,0)(1,1.72)
\psline(1,1.72)(-2,0)\psline(1,1.72)(1,-1.72)\psline(-1,1.72)(-1,-1.72)
\psline(-1,1.72)(2,0)\psline(-2,0)(1,-1.72)\psline(-2,0)(1,1.72)
\psline(-1,-1.72)(2,0)\psline(-1,-1.72)(1,1.72)\psline(2,0)(1,1.72)\psline(2,0)(-2,0)
\pscircle[fillcolor=lightgray, fillstyle=solid, linewidth=1pt](1,1.72){.3}
\pscircle[fillcolor=lightgray, fillstyle=solid, linewidth=1pt](-1,1.72){.3}
\pscircle[fillcolor=lightgray, fillstyle=solid, linewidth=1pt](-2,0){.3}
\pscircle[fillcolor=lightgray, fillstyle=solid, linewidth=1pt](-1,-1.72){.3}
\pscircle[fillcolor=lightgray, fillstyle=solid, linewidth=1pt](1,-1.72){.3}
\pscircle[fillcolor=lightgray, fillstyle=solid, linewidth=1pt](2,0){.3}
\pscircle[fillcolor=white, fillstyle=solid, linecolor=white, linewidth=0pt](0.7,1.2){.15}
\pscircle[fillcolor=white, fillstyle=solid, linecolor=white, linewidth=0pt](1,1.1){.15}
\pscircle[fillcolor=white, fillstyle=solid, linecolor=white, linewidth=0pt](0,1.72){.15}
\pscircle[fillcolor=white, fillstyle=solid, linecolor=white, linewidth=0pt](0.43,1.42){.15}
\pscircle[fillcolor=white, fillstyle=solid, linecolor=white, linewidth=0pt](1.42,1.1){.15}
\pscircle[fillcolor=white, fillstyle=solid, linecolor=white, linewidth=0pt](-0.7,1.2){.15}
\pscircle[fillcolor=white, fillstyle=solid, linecolor=white, linewidth=0pt](-1,1.1){.15}
\pscircle[fillcolor=white, fillstyle=solid, linecolor=white, linewidth=0pt](-0.43,1.42){.15}
\pscircle[fillcolor=white, fillstyle=solid, linecolor=white, linewidth=0pt](-1.42,1.1){.15}
\pscircle[fillcolor=white, fillstyle=solid, linecolor=white, linewidth=0pt](-0.7,-1.2){.15}
\pscircle[fillcolor=white, fillstyle=solid, linecolor=white, linewidth=0pt](-1,-1.1){.15}
\pscircle[fillcolor=white, fillstyle=solid, linecolor=white, linewidth=0pt](-0.43,-1.42){.15}
\pscircle[fillcolor=white, fillstyle=solid, linecolor=white, linewidth=0pt](-1.42,-1.1){.15}
\pscircle[fillcolor=white, fillstyle=solid, linecolor=white, linewidth=0pt](0.7,-1.2){.15}
\pscircle[fillcolor=white, fillstyle=solid, linecolor=white, linewidth=0pt](1,-1.1){.15}
\pscircle[fillcolor=white, fillstyle=solid, linecolor=white, linewidth=0pt](0,-1.72){.15}
\pscircle[fillcolor=white, fillstyle=solid, linecolor=white, linewidth=0pt](0.43,-1.42){.15}
\pscircle[fillcolor=white, fillstyle=solid, linecolor=white, linewidth=0pt](1.42,-1.1){.15}
\pscircle[fillcolor=white, fillstyle=solid, linecolor=white, linewidth=0pt](1.38,0){.15}
\pscircle[fillcolor=white, fillstyle=solid, linecolor=white, linewidth=0pt](1.42,0.3){.15}
\pscircle[fillcolor=white, fillstyle=solid, linecolor=white, linewidth=0pt](1.42,-0.3){.15}
\pscircle[fillcolor=white, fillstyle=solid, linecolor=white, linewidth=0pt](-1.38,0){.15}
\pscircle[fillcolor=white, fillstyle=solid, linecolor=white, linewidth=0pt](-1.42,0.3){.15}
\pscircle[fillcolor=white, fillstyle=solid, linecolor=white, linewidth=0pt](-1.42,-0.3){.15}
\rput(1,1.72){{$2$}}\rput(-1,1.72){{$5$}}\rput(-2,0){{$3$}}
\rput(-1,-1.72){{$1$}}\rput(1,-1.72){{$6$}}\rput(2,0){{$4$}}
\rput(0,1.72){{\footnotesize{7}}}\rput(0.43,1.42){{\footnotesize{6}}}
\rput(0.7,1.2){{\footnotesize{5}}}
\rput(1,1.1){{\footnotesize{4}}}\rput(1.42,1.1){{\footnotesize{3}}}
\rput(-0.43,1.42){{\footnotesize{1}}}\rput(-0.7,1.2){{\footnotesize{2}}}
\rput(-1,1.1){{\footnotesize{3}}}\rput(-1.42,1.1){{\footnotesize{4}}}
\rput(0,-1.72){{\footnotesize{7}}}\rput(-0.43,-1.42){{\footnotesize{6}}}
\rput(-0.7,-1.2){{\footnotesize{5}}}\rput(-1,-1.1){{\footnotesize{3}}}
\rput(-1.42,-1.1){{\footnotesize{2}}}
\rput(0.43,-1.42){{\footnotesize{1}}}\rput(0.7,-1.2){{\footnotesize{2}}}
\rput(1,-1.1){{\footnotesize{4}}}\rput(1.42,-1.1){{\footnotesize{5}}}
\rput(1.42,0.3){{\footnotesize{1}}}\rput(1.38,0){{\footnotesize{7}}}
\rput(1.42,-0.3){{\footnotesize{6}}}\rput(-1.42,0.3){{\footnotesize{6}}}
\rput(-1.38,0){{\footnotesize{7}}}\rput(-1.42,-0.3){{\footnotesize{1}}}
\rput(1.35,2.32){{$\{1\}^{c}$}}\rput(2.7,0){{$\{2\}^{c}$}}
\rput(1.35,-2.32){{$\{3\}^{c}$}}\rput(-1.35,-2.32){{$\{4\}^{c}$}}
\rput(-2.7,0){{$\{5\}^{c}$}}\rput(-1.5,2.32){{$\{6\}^{c}$}} \rput(0,-2.9){{$(b)$}}
\end{pspicture}
$$
\caption{Proper total colorings distinguishing adjacent vertices
by sums of the complete graph $(a)$ $K_5$ and $(b)$ $K_6$.} \label{k6fig}
\end{figure}
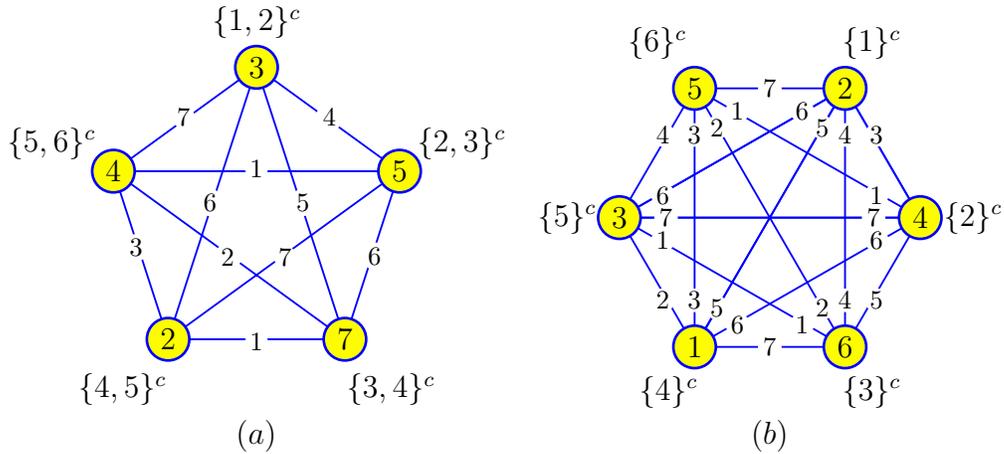

M. pil\'{s}niak and M. Wo\'{z}niak~\cite{PW} first introduced that a proper total coloring
of $\Gamma$ is \emph{a proper total colorings distinguishing
adjacent vertices by sums} if for a vertex $v \in V(\Gamma)$,
the total sum of colors of the edges incident to $v$ and the color
of $v$, denoted by $f(v)$, are distinct for adjacent vertices.
An example of a proper total $7$ coloring distinguishing adjacent
vertices by sums of the complete graph $K_6$ is given Figure~\ref{k6fig}.
The smallest number of color $k$
such that $\Gamma$ admits a proper total $k$ colorings distinguishing adjacent vertices by sums is
called \emph{adjacent vertex distinguishing index by sum}, denoted by $\tndi (\Gamma)$.
They also find that the adjacent vertex distinguishing index by sums of some special graphs
including paths, cycles, stars, complete and complete bipartite graphs.
They made a well-believed conjecture that

\begin{conj}\cite{PW} \label{conjecture} Let $G$ be a graph with the maximum degree $\Delta$.
Then,
$$\tndi (G)\le \Delta +3.$$
\end{conj}

Consequently, they showed that the above graphs hold Conjecture~\ref{conjecture} and they also
proved that the conjecture holds for regular bipartite graphs, cubic graphs, graphs with $\Delta \le 3$.
Li, Liu and Wang~\cite{LLW} proved that Conjecture~\ref{conjecture} holds for $K_4$-minor free graphs.

It is obvious that $ \Delta +1 \le \tndi (G)$. Thus if the conjecture
is true, then we can divide graphs into three groups;
$\tndi (G)= \Delta +1$, $\Delta +2$ and $\Delta +3$,
we called them $\tndi$ class I, II and III, respectively.

Once some types of graphs were studied, it is natural to consider
their products for the next step.
Chen, Zhang and Sun~\cite{CZS} found the adjacent vertex
distinguishing total chromatic numbers
of $P_m \times P_n$, $P_m \times C_n$ and $C_m \times C_n$.

The aim of the present article is two-fold. First we find that
the wheel graphs are $\tndi$ class I except $W_3$ which is $K_4$.
Second, we find that the product graphs $P_m \times P_n$, $C_m \times P_n$,
$S_m \times P_n$, $W_m \times P_n$ and $K_m \times P_n$, $S_m \times C_n$ and
$C_m \times C_n$ are $\tndi$ class II except $P_3 \times P_3$, $S_m \times P_3$
and $W_m \times P_3$ which are $\tndi$ class I.
These results not only find the adjacent vertex distinguishing indices by sums
but also strongly support that Conjecture~\ref{conjecture}
holds for these product graphs.

The outline of this paper is as follows.
We first provide some preliminary definitions and results
in section~\ref{prelim}. In section~\ref{other},
we find the adjacent vertex distinguishing indices by sums of
the star graphs and wheel graphs.
In section~\ref{prod}, we investigate the adjacent vertex distinguishing indices by sum of
the products of two graphs.
At last, we make a conclusive remark in section~\ref{endcom}.

\section{Preliminaries} \label{prelim}

For a graph $\Gamma$, the maximal and minimal degree of $\Gamma$
are denoted by $\Delta$ and $\delta$ respectively.
The following theorems will be used in proofs of theorems in section~\ref{other} and~\ref{prod}.
For general terminology in graph theory, we refer the reader
to~\cite{BM}.

\begin{theorem}(\cite[Observation 3]{PW}) If a graph $\Gamma$
contains two adjacent vertices $x$, $y$ such that $deg(x)=deg(y)=\Delta$, then
$\tndi (\Gamma) \ge \Delta+2$. \label{obs1}
\end{theorem}


\begin{theorem} (\cite[Proposition 4, 5, 6, 7, 8, 9 $\&$ 10]{PW}) \label{tndioldthm}
\begin{enumerate}
\item Let $P_n$ be a path of the size $n\ge 2$. Then
$$ \tndi (P_n) = \left\{ \begin{array}{cl}
\Delta+1 & {\rm{if}}~ n=3, \\
\Delta+2 & {\rm{if}}~ n=2~{\rm{or}}~ n\ge 4. \end{array}\right.$$
\item Let $C_n$ be a cycle of the size $n\ge 3$. Then
$$ \tndi  (C_n) = \left\{ \begin{array}{cl}
\Delta+3& \rm{if}~ \textit{n}=3, \\
\Delta+2& \rm{if}~\textit{n}\ge 4. \end{array}\right.$$
\item Let $S_n$ be a star graph of size $n\ge 2$.
Then
$$ \tndi  (S_n) = \Delta+1.$$
\item Let $K_n$ be a complete graph of the size $n\ge 2$. Then
$$ \tndi  (K_n) = \left\{ \begin{array}{cl}
\Delta+2 & \rm{if} ~ \textit{n}~\rm{is}~\rm{even}, \\
\Delta+3 & \rm{if} ~ \textit{n}~\rm{is}~\rm{odd}. \end{array}\right.$$
\item Let $K_{p,q}$ be a complete bipartite graph. Then
$$ \tndi  (K_{m,n}) = \left\{ \begin{array}{cl}
\Delta+1& {\rm{if}} ~ n \neq m, \\
\Delta+2 & {\rm{if}}~ n=m. \end{array}\right.$$
\item Let $\Gamma$ be a regular bipartite graph. Then $\tndi (\Gamma) = \Delta+2$.
\item Let $T$ be a tree of order $n$. Then $\Delta+1 \le \tndi (T) \le \Delta+2$.
Furthermore, if there exist two adjacent vertices $x$, $y$ such that $deg(x)=deg(y)=\Delta$, then
$\tndi  (T) = \Delta+2$, and otherwise, we have $\tndi  (T) = \Delta+1$.
\end{enumerate}
\end{theorem}


Throughout the article, we will often use $\equiv$ (mod $n$).
Since this modulo will be used for colorings,
we use the standard complete residue system which is
$\{ 1, 2, \ldots, n-1, n\}$ for modulo $n$, unless stated differently.

\section{Some other special graphs} \label{other}

In this section, we will discuss the adjacent vertex distinguishing indices by sums
of the star graph and the wheel graphs. Although the adjacent vertex distinguishing indices by sums
of the star graphs are already known in Theorem~\ref{tndioldthm} (3), we will find
two proper total colorings distinguishing adjacent vertices by sums for
these two graphs because they will be used in section~\ref{prod}.

\begin{theorem}  \label{snthm}
Let $S_n$ be a star graph of size $n\ge 2$.
Then
$$ \tndi (S_n) = \Delta+1.$$
\begin{proof}
A proper total coloring $c$ distinguishing adjacent vertices by sums
provided in Theorem~\ref{tndioldthm} is as follows ;
the edge of $S_n$ is colored by $1, 2, \ldots, \Delta$,
the vertex in the middle is colored by $\Delta+1$,
the vertex incident to the edge colored by $1$ is colored by $2$
and all remaining vertices are colored by $1$.

For the second proper total coloring $c'$ distinguishing adjacent vertices by sums of $S_n$,
we consider $c'\equiv c+1$ modulo $\Delta+1$. One can easily check this is
also a proper total coloring distinguishing adjacent vertices by sums of $S_n$.
\end{proof}
\end{theorem}

The wheel graphs are highly well structured, for example, they are planar
and any maximal planar graph, other than $K_4$, contains a subgraph isomorphic to
either $W_5$ or $W_6$. For $n \ge 3$, the wheel graph $W_n$ contains $n^2-3n+3$ cycles
and at least one of them is a Hamiltonian cycle. Their (vertex) chromatic number are already known as

$$ \chi (W_m) = \left\{ \begin{array}{cl}
3& \rm{if} ~\textit{m}~\rm{is}~\rm{odd}, \\
4 & \rm{if} ~\textit{m}~\rm{is}~\rm{even}. \end{array}\right.$$

The following theorem finds the proper total coloring distinguishing adjacent vertices by sums of the wheel graph.

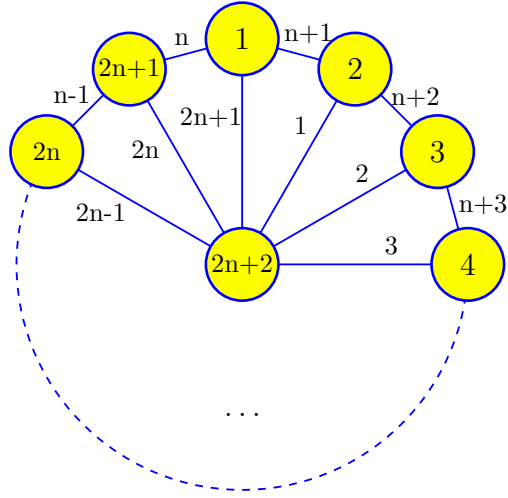
\begin{figure}
$$
\begin{pspicture}[shift=-2](-3.6,-3)(3.6,3.6)
\psline(0,0)(3;0) \psline(0,0)(3;30) \psline(0,0)(3;60) \psline(0,0)(3;90)
\psline(0,0)(3;120) \psline(0,0)(3;150)
\psline(3;0)(3;30)(3;60)(3;90)(3;120)(3;150)
\psarc[linestyle=dashed](0,0){3}{150}{0}
\pscircle[fillcolor=lightgray, fillstyle=solid, linewidth=1pt](0,0){.5}
\pscircle[fillcolor=lightgray, fillstyle=solid, linewidth=1pt](3;0){.5}
\pscircle[fillcolor=lightgray, fillstyle=solid, linewidth=1pt](3;30){.5}
\pscircle[fillcolor=lightgray, fillstyle=solid, linewidth=1pt](3;60){.5}
\pscircle[fillcolor=lightgray, fillstyle=solid, linewidth=1pt](3;90){.5}
\pscircle[fillcolor=lightgray, fillstyle=solid, linewidth=1pt](3;120){.5}
\pscircle[fillcolor=lightgray, fillstyle=solid, linewidth=1pt](3;150){.5}
\rput(3;90){{$1$}} \rput(3;60){{$2$}} \rput(3;30){{$3$}} \rput(3;0){{$4$}}
\rput(3;120){{\footnotesize{2n+1}}} \rput(3;150){{\footnotesize{2n}}}
\rput(2;67){{\footnotesize{1}}} \rput(2;37){{\footnotesize{2}}} \rput(2;7){{\footnotesize{3}}}
\rput(2;102){{\footnotesize{2n+1}}} \rput(2;130){{\footnotesize{2n}}}
\rput(2;160){{\footnotesize{2n-1}}}
\rput(3.2;135){{\footnotesize{n-1}}}\rput(3.1;105){{\footnotesize{n}}}
\rput(3.2;74){{\footnotesize{n+1}}}\rput(3.2;44){{\footnotesize{n+2}}}
\rput(3.3;14){{\footnotesize{n+3}}}
\rput(0,0){{\footnotesize{2n+2}}} \rput(0,-2){{$\ldots$}}
\end{pspicture}
$$
\caption{A proper total coloring distinguishing adjacent vertices
by sums of the wheel graph $W_{2n+1}$, $n \ge 2$.} \label{W4odd}
\end{figure}

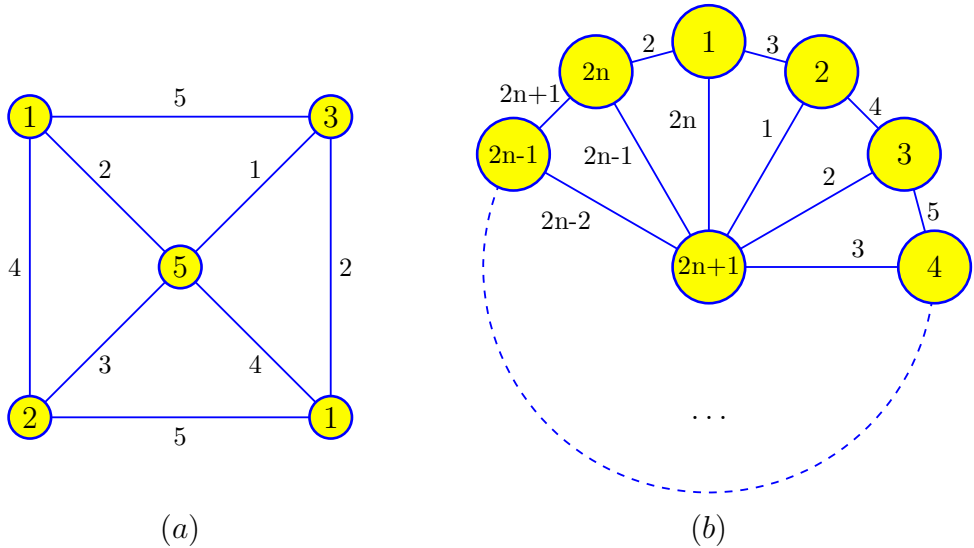
\begin{figure}
$$
\begin{pspicture}[shift=-2](-2.6,-2.6)(2.6,2.6)
\psline(-2,2)(2,-2)
\psline(2,2)(-2,-2)
\psline(2,2)(2,-2)(-2,-2)(-2,2)(2,2)
\pscircle[fillcolor=lightgray, fillstyle=solid, linewidth=1pt](2,2){.3}
\pscircle[fillcolor=lightgray, fillstyle=solid, linewidth=1pt](2,-2){.3}
\pscircle[fillcolor=lightgray, fillstyle=solid, linewidth=1pt](-2,2){.3}
\pscircle[fillcolor=lightgray, fillstyle=solid, linewidth=1pt](-2,-2){.3}
\pscircle[fillcolor=lightgray, fillstyle=solid, linewidth=1pt](0,0){.3}
\rput(0,0){{$5$}} \rput(-2,2){{$1$}} \rput(2,2){{$3$}}
\rput(2,-2){{$1$}} \rput(-2,-2){{$2$}}
\rput(2.2,0){{\footnotesize{2}}} \rput(-2.2,0){{\footnotesize{4}}}
\rput(0,2.25){{\footnotesize{5}}} \rput(0,-2.25){{\footnotesize{5}}}
\rput(1,1.3){{\footnotesize{1}}} \rput(-1,1.3){{\footnotesize{2}}}
\rput(1,-1.3){{\footnotesize{4}}} \rput(-1,-1.3){{\footnotesize{3}}}
\rput(0,-3.5){{$(a)$}}
\end{pspicture}
\quad\quad
\begin{pspicture}[shift=-3](-3.6,-3.6)(3.6,3.6)
\psline(0,0)(3;0) \psline(0,0)(3;30) \psline(0,0)(3;60) \psline(0,0)(3;90)
\psline(0,0)(3;120) \psline(0,0)(3;150)
\psline(3;0)(3;30)(3;60)(3;90)(3;120)(3;150)
\psarc[linestyle=dashed](0,0){3}{150}{0}
\pscircle[fillcolor=lightgray, fillstyle=solid, linewidth=1pt](0,0){.5}
\pscircle[fillcolor=lightgray, fillstyle=solid, linewidth=1pt](3;0){.5}
\pscircle[fillcolor=lightgray, fillstyle=solid, linewidth=1pt](3;30){.5}
\pscircle[fillcolor=lightgray, fillstyle=solid, linewidth=1pt](3;60){.5}
\pscircle[fillcolor=lightgray, fillstyle=solid, linewidth=1pt](3;90){.5}
\pscircle[fillcolor=lightgray, fillstyle=solid, linewidth=1pt](3;120){.5}
\pscircle[fillcolor=lightgray, fillstyle=solid, linewidth=1pt](3;150){.5}
\rput(3;90){{$1$}} \rput(3;60){{$2$}} \rput(3;30){{$3$}} \rput(3;0){{$4$}}
\rput(3;120){{\footnotesize{2n}}} \rput(3;150){{\footnotesize{2n-1}}}
\rput(2;67){{\footnotesize{1}}} \rput(2;37){{\footnotesize{2}}} \rput(2;7){{\footnotesize{3}}}
\rput(2;100){{\footnotesize{2n}}} \rput(2;132){{\footnotesize{2n-1}}}
\rput(2;162){{\footnotesize{2n-2}}}
\rput(3.3;136){{\footnotesize{2n+1}}}\rput(3.08;105){{\footnotesize{2}}}
\rput(3.08;74){{\footnotesize{3}}}\rput(3.08;44){{\footnotesize{4}}}
\rput(3.08;14){{\footnotesize{5}}}
\rput(0,0){{\footnotesize{2n+1}}} \rput(0,-2){{$\ldots$}}
\rput(0,-3.5){{$(b)$}}
\end{pspicture}
$$
\caption{A proper total coloring distinguishing adjacent vertices by
sums of the wheel graph $(a)$ $W_4$ and $(b)$ $W_{2n}$, $n \ge 3$.} \label{W4even}
\end{figure}

\begin{theorem} \label{wnthm}
Let $W_m$ be a wheel graph of size $m\ge 3$.
Then
$$ \tndi (W_m) = \left\{ \begin{array}{cl}
\Delta+1 & {\rm{for}}~ m\ge 4, \\
\Delta+2 & {\rm{for}}~ m=3. \end{array}\right.$$
\begin{proof}
The smallest $m$ for which the wheel graph makes a natural geometric shape is $3$.
But, because $W_3 \equiv K_4$ and $\tndi (K_4)$ was given in Theorem~\ref{tndioldthm} (4).
Now we assume $m\ge 4$. We prove the theorem dividing cases depending on the parity of $m$.

First we assume $m$ is odd, so write $m=2n+1$ for some $n\ge 2$.
We color the center vertex of the wheel graph by color $2n+2$
and the adjacent edges by $1, 2, \ldots, 2n+1(=\Delta)$ clockwisely.
Next we use color $2$ for a vertex connected to the center vertex by the edge colored by $1$,
and for all remaining vertices, we similarly color by $3, \ldots, 2n+1, 1$ clockwisely.
Last we will color the edges in the circumference of the wheel graph.
The edge between vertices colored by $1$ and $2$ will be colored by $n+1$,
and we increase the color by $1$ up to mod$(2n+1)$ clockwisely.
Therefore, the wheel graph of an odd size has $\tndi (W_m) \le m+1$ as illustrated in Figure~\ref{W4odd}.
By combining the fact $m+1=\Delta +1 \le \tndi (W_m)$, we find $\tndi (W_m) = m+1$.

Next, we assume $m$ is even. One may find $\tndi (W_4) = 5$ as shown in Figure~\ref{W4even} $(a)$.
Now we assume that $m=2n$ is bigger than or equal to $6$.
Similar to the $m=2n+1$ we color the vertices and edges
as illustrated in Figure~\ref{W4even} $(b)$.
By the same reason, we find $\tndi (W_n) =\Delta+1$.
It completes the proof of the theorem.

Let $c$ be the coloring of $W_m$ discussed above.
For the second proper total coloring $c'$ distinguishing adjacent vertices by sums of $W_m$,
we consider $c'\equiv c+1$ modulo $\Delta+1$ ($\Delta +2$, resp) for $m\ge 4$ ($m=3$, respectively).
\end{proof}
\end{theorem}

Although we are not discussing the hypercubes,
one may find it is regular bipartite. So
by using Theorem~\ref{tndioldthm} (6), it is $\tndi$ class II.
The ladder graph can be considered as a product of two graphs and will be handled in the following section.

\section{Products of two graphs} \label{prod}

Throughout this section, the product means the graph product unless stated differently.
A formal definition of the product is as follows.
Let $G_1=(V(G_1), E(G_1))$, $G_2=(V(G_2), E(G_2))$ be two graphs.
The \emph{product of $G_1$ and $G_2$} which is denoted by $G_1 \times G_2$ is
consists of the vertex set
$$V(G_1 \times G_2)= V(G_1) \times V(G_2)$$
and the edge set
\begin{align*}
E(G_1 \times G_2) = \{ ((u_1, u_2), (v_1, v_2)) | &\left[u_1 = v_1 ~{\rm{and}}~ (u_2, v_2) \in E(G_2)\right]
~{\rm{or}}\\
&\left[(u_1, v_1) \in E(G_1) ~{\rm{and}}~ u_2 = v_2 \right] \}.
\end{align*}

Throughout the section, we use the following notations;
let $V(G_2) = \{ v_1, v_2, \ldots, v_l \}$.
$G_1 \times \{v_i \}$ is \emph{the $i$-th copy of $G_1$} denoted by $(G_1 \times G_2)_i$,
and the set of the edges between $(G_1 \times G_2)_i$ and $(G_1 \times G_2)_{i+1}$
is denoted by $h((G_1 \times G_2)_i)$ where $i=1, 2, \ldots, l$.

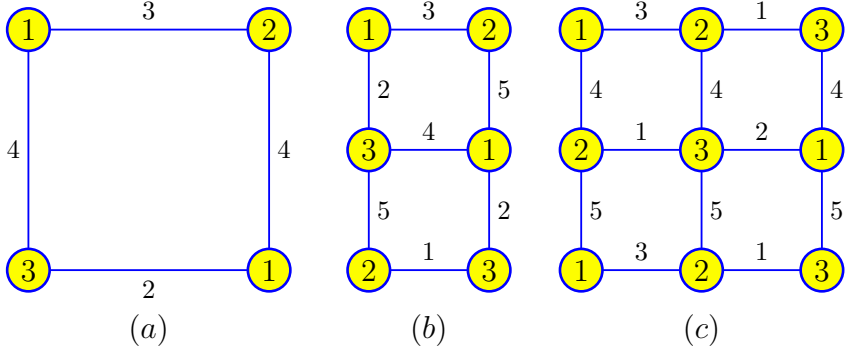
\begin{figure}
$$
\begin{pspicture}[shift=-2](-1.8,-2.7)(1.8,1.8)
\psline(1.6,1.6)(1.6,-1.6)(-1.6,-1.6)(-1.6,1.6)(1.6,1.6)
\pscircle[fillcolor=lightgray, fillstyle=solid, linewidth=1pt](1.6,1.6){.3}
\pscircle[fillcolor=lightgray, fillstyle=solid, linewidth=1pt](1.6,-1.6){.3}
\pscircle[fillcolor=lightgray, fillstyle=solid, linewidth=1pt](-1.6,1.6){.3}
\pscircle[fillcolor=lightgray, fillstyle=solid, linewidth=1pt](-1.6,-1.6){.3}
\rput(-1.6,1.6){{$1$}} \rput(1.6,1.6){{$2$}} \rput(1.6,-1.6){{$1$}} \rput(-1.6,-1.6){{$3$}}
\rput(1.8,0){{\footnotesize{4}}} \rput(-1.8,0){{\footnotesize{4}}}
\rput(0,1.85){{\footnotesize{3}}} \rput(0,-1.85){{\footnotesize{2}}}
\rput(0,-2.4){{$(a)$}}
\end{pspicture}
\quad\quad
\begin{pspicture}[shift=-2](-.3,-2.7)(1.8,1.8)
\psline(1.6,1.6)(1.6,-1.6)(0,-1.6)(0,1.6)(1.6,1.6)
\psline(0,0)(1.6,0)
\pscircle[fillcolor=lightgray, fillstyle=solid, linewidth=1pt](1.6,1.6){.3}
\pscircle[fillcolor=lightgray, fillstyle=solid, linewidth=1pt](1.6,-1.6){.3}
\pscircle[fillcolor=lightgray, fillstyle=solid, linewidth=1pt](0,0){.3}
\pscircle[fillcolor=lightgray, fillstyle=solid, linewidth=1pt](1.6,0){.3}
\pscircle[fillcolor=lightgray, fillstyle=solid, linewidth=1pt](0,1.6){.3}
\pscircle[fillcolor=lightgray, fillstyle=solid, linewidth=1pt](0,-1.6){.3}
\rput(1.6,1.6){{$2$}} \rput(1.6,-1.6){{$3$}} 
\rput(1.6,0){{$1$}} \rput(0,1.6){{$1$}} \rput(0,-1.6){{$2$}}
\rput(0,0){{$3$}}
\rput(1.8,0.8){{\footnotesize{5}}}
\rput(0.2,0.8){{\footnotesize{2}}}\rput(1.8,-0.8){{\footnotesize{2}}}
\rput(0.2,-0.8){{\footnotesize{5}}}
\rput(0.8,1.85){{\footnotesize{3}}}\rput(0.8,-1.35){{\footnotesize{1}}}
\rput(0.8,0.25){{\footnotesize{4}}}
\rput(.8,-2.4){{$(b)$}}
\end{pspicture}
\quad\quad
\begin{pspicture}[shift=-2](-1.8,-2.7)(1.8,1.8)
\psline(1.6,1.6)(1.6,-1.6)(-1.6,-1.6)(-1.6,1.6)(1.6,1.6)
\psline(-1.6,0)(1.6,0)
\psline(0,1.6)(0,-1.6)
\pscircle[fillcolor=lightgray, fillstyle=solid, linewidth=1pt](1.6,1.6){.3}
\pscircle[fillcolor=lightgray, fillstyle=solid, linewidth=1pt](1.6,-1.6){.3}
\pscircle[fillcolor=lightgray, fillstyle=solid, linewidth=1pt](-1.6,1.6){.3}
\pscircle[fillcolor=lightgray, fillstyle=solid, linewidth=1pt](-1.6,-1.6){.3}
\pscircle[fillcolor=lightgray, fillstyle=solid, linewidth=1pt](0,0){.3}
\pscircle[fillcolor=lightgray, fillstyle=solid, linewidth=1pt](1.6,0){.3}
\pscircle[fillcolor=lightgray, fillstyle=solid, linewidth=1pt](0,1.6){.3}
\pscircle[fillcolor=lightgray, fillstyle=solid, linewidth=1pt](0,-1.6){.3}
\pscircle[fillcolor=lightgray, fillstyle=solid, linewidth=1pt](-1.6,0){.3}
\rput(-1.6,1.6){{$1$}} \rput(1.6,1.6){{$3$}} \rput(1.6,-1.6){{$3$}} \rput(-1.6,-1.6){{$1$}}
\rput(1.6,0){{$1$}} \rput(-1.6,0){{$2$}}\rput(0,1.6){{$2$}} \rput(0,-1.6){{$2$}}\rput(0,0){{$3$}}
\rput(1.8,0.8){{\footnotesize{4}}}\rput(-1.4,0.8){{\footnotesize{4}}}
\rput(0.2,0.8){{\footnotesize{4}}}\rput(1.8,-0.8){{\footnotesize{5}}}
\rput(-1.4,-0.8){{\footnotesize{5}}}\rput(0.2,-0.8){{\footnotesize{5}}}
\rput(0.8,1.85){{\footnotesize{1}}}\rput(0.8,-1.35){{\footnotesize{1}}}
\rput(0.8,0.25){{\footnotesize{2}}}\rput(-0.8,1.85){{\footnotesize{3}}}
\rput(-0.8,-1.35){{\footnotesize{3}}}\rput(-0.8,0.25){{\footnotesize{1}}}
\rput(0,-2.4){{$(c)$}}
\end{pspicture}
$$
\caption{A proper total coloring distinguishing adjacent vertices by
sums of the product graph $(a)$ $P_2 \times P_2$, $(b)$ $P_2 \times P_3$ and $(c)$ $P_3 \times P_3$.} \label{p2pmfig}
\end{figure}

\begin{theorem}
Let $P_m, P_n$ be a path of order $m, n \ge 2$. Then
$$ \tndi (P_m \times P_n) = \left\{ \begin{array}{cl}
\Delta +1 & {\rm{for}} ~m=n=3, \\
\Delta +2 & {\rm{Otherwise}}. \end{array}\right.$$
\label{pnpm}
\begin{proof}
First, if $m=n=2$ or $3$, then maximum degree $\Delta$ is $2$ or $4$, respectively.
Proper total colorings distinguishing adjacent vertices by sums
of $P_2 \times P_2$ and $P_3 \times P_3$ are depicted in Figure~\ref{p2pmfig}.
For the case $P_2\times P_3$, then maximum degree $\Delta$ is $3$ and two
adjacent vertices have the maximum degree $3$, by Theorem~\ref{obs1} and 
the coloring in Figure~\ref{p2pmfig} $(b)$, we find that $\tndi (P_2 \times P_3) = 5$.
Let us remark that $P_3 \times P_3$ is the unique case that there is only one
vertex of the maximum degree and 
$\tndi(P_3 \times P_3)=5=\Delta+1$ which is the only exception
that $P_m \times P_n$ is not $\tndi$ class II.

Second, if just one of $m, n$ is $2$, then without loss of generality, we may assume $n=2$ and $m\ge 3$.
The maximum degree $\Delta$ of $P_m \times P_2$ is $3$.
Color the vertices of $(P_m \times P_2)_1$ alternatively with colors $1, 3, 1, 3, \ldots$,
and the edges with $2, 4, 2, 4, \ldots$.
For $(P_m \times P_2)_2$, we first color the edges with $4, 2, 4, 2, \ldots$, 
and the vertices $2, 1, 3, 1, \ldots$, $2(4)$ if the last color of edges was $4(2, reps.)$.
as two different rightmost figures as depicted in Figure~\ref{pnpmfig}.
Color the edges in $h((P_m \times P_2)_1)$ by $5$ then this is a proper
total $5$ coloring distinguishing adjacent vertices by sums.
We get $\tndi(P_m \times P_2) \le 5$ but Theorem~\ref{obs1} implies
that $5= \Delta + 2 \le \tndi(P_m \times P_2)$. Thus, $\tndi(P_m \times P_2) = \Delta +2$.
This coloring will be used for the remaining products $P_m \times P_n$, let us denote it by $\phi$.

\begin{figure}
$$
\begin{pspicture}[shift=-2](-3.5,-3)(6.2,2)
\psline(-3.2,1.6)(2.4,1.6) \psline(1.6,1.6)(1.6,-2.4)
\psline(-3.2,-1.6)(2.4,-1.6) \psline(-3.2,1.6)(-3.2,-2.4)
\psline(-3.2,0)(2.4,0) \psline(-1.6,1.6)(-1.6,-2.4)
\psline(0,1.6)(0,-2.4) \psline(3.4,1.6)(4.2,1.6)
\psline(3.4,0)(4.2,0) \psline(3.4,-1.6)(4.2,-1.6)
\psline(4.2,1.6)(4.2,-2.4) \psline(5.2,1.6)(6,1.6)
\psline(5.2,0)(6,0) \psline(5.2,-1.6)(6,-1.6)
\psline(6,1.6)(6,-2.4)
\pscircle[fillcolor=lightgray, fillstyle=solid, linewidth=1pt](-3.2,1.6){.3}
\pscircle[fillcolor=lightgray, fillstyle=solid, linewidth=1pt](-1.6,1.6){.3}
\pscircle[fillcolor=lightgray, fillstyle=solid, linewidth=1pt](0,1.6){.3}
\pscircle[fillcolor=lightgray, fillstyle=solid, linewidth=1pt](1.6,1.6){.3}
\pscircle[fillcolor=lightgray, fillstyle=solid, linewidth=1pt](-3.2,0){.3}
\pscircle[fillcolor=lightgray, fillstyle=solid, linewidth=1pt](-1.6,0){.3}
\pscircle[fillcolor=lightgray, fillstyle=solid, linewidth=1pt](0,0){.3}
\pscircle[fillcolor=lightgray, fillstyle=solid, linewidth=1pt](1.6,0){.3}
\pscircle[fillcolor=lightgray, fillstyle=solid, linewidth=1pt](-3.2,-1.6){.3}
\pscircle[fillcolor=lightgray, fillstyle=solid, linewidth=1pt](-1.6,-1.6){.3}
\pscircle[fillcolor=lightgray, fillstyle=solid, linewidth=1pt](0,-1.6){.3}
\pscircle[fillcolor=lightgray, fillstyle=solid, linewidth=1pt](1.6,-1.6){.3}
\pscircle[fillcolor=lightgray, fillstyle=solid, linewidth=1pt](4.2,1.6){.3}
\pscircle[fillcolor=lightgray, fillstyle=solid, linewidth=1pt](4.2,0){.3}
\pscircle[fillcolor=lightgray, fillstyle=solid, linewidth=1pt](4.2,-1.6){.3}
\pscircle[fillcolor=lightgray, fillstyle=solid, linewidth=1pt](6,1.6){.3}
\pscircle[fillcolor=lightgray, fillstyle=solid, linewidth=1pt](6,0){.3}
\pscircle[fillcolor=lightgray, fillstyle=solid, linewidth=1pt](6,-1.6){.3}
\rput(-3.2,1.6){{$1$}} \rput(-1.6,1.6){{$3$}} \rput(0,1.6){{$1$}} \rput(1.6,1.6){{$3$}}
\rput(-3.2,0){{$2$}} \rput(-1.6,0){{$1$}} \rput(0,0){{$3$}} \rput(1.6,0){{$1$}}
\rput(-3.2,-1.6){{$1$}}\rput(-1.6,-1.6){{$3$}}\rput(0,-1.6){{$1$}}\rput(1.6,-1.6){{$3$}}
\rput(4.2,1.6){{$3$}}\rput(4.2,0){{$2$}}\rput(4.2,-1.6){{$3$}}  
\rput(6,1.6){{$1$}}\rput(6,0){{$4$}}\rput(6,-1.6){{$1$}} 
\rput(1.76,0.8){{\footnotesize{5}}}\rput(-1.44,0.8){{\footnotesize{5}}}
\rput(0.16,0.8){{\footnotesize{5}}}\rput(-3.04,0.8){{\footnotesize{5}}}
\rput(4.36,0.8){{\footnotesize{5}}}\rput(6.16,0.8){{\footnotesize{5}}}
\rput(1.76,-0.8){{\footnotesize{6}}}\rput(-1.44,-0.8){{\footnotesize{6}}}
\rput(0.16,-0.8){{\footnotesize{6}}}\rput(-3.04,-0.8){{\footnotesize{6}}}
\rput(4.36,-0.8){{\footnotesize{6}}}\rput(6.16,-0.8){{\footnotesize{6}}}
\rput(1.76,-2.2){{\footnotesize{5}}}\rput(-1.44,-2.2){{\footnotesize{5}}}
\rput(0.16,-2.2){{\footnotesize{5}}}\rput(-3.04,-2.2){{\footnotesize{5}}}
\rput(4.36,-2.2){{\footnotesize{5}}}\rput(6.16,-2.2){{\footnotesize{5}}}
\rput(0.8,1.85){{\footnotesize{2}}}\rput(0.8,-1.35){{\footnotesize{2}}}
\rput(0.8,0.25){{\footnotesize{4}}}
\rput(-0.8,1.85){{\footnotesize{4}}}\rput(-0.8,-1.35){{\footnotesize{4}}}
\rput(-0.8,0.25){{\footnotesize{2}}}
\rput(-2.4,1.85){{\footnotesize{2}}}\rput(-2.4,-1.35){{\footnotesize{2}}}
\rput(-2.4,0.25){{\footnotesize{4}}}
\rput(2.2,1.85){{\footnotesize{4}}}\rput(2.2,-1.35){{\footnotesize{4}}}
\rput(2.2,0.25){{\footnotesize{2}}}
\rput(3.6,1.85){{\footnotesize{2}}}\rput(3.6,-1.35){{\footnotesize{2}}}
\rput(3.6,0.25){{\footnotesize{4}}}
\rput(5.4,1.85){{\footnotesize{4}}}\rput(5.4,-1.35){{\footnotesize{4}}}
\rput(5.4,0.25){{\footnotesize{2}}}
\rput(2.8,1.6){{$\cdots$}} \rput(2.8,0){{$\cdots$}} \rput(2.8,-1.6){{$\cdots$}}
\rput(-3.2,-2.8){{$\vdots$}} \rput(-1.6,-2.8){{$\vdots$}}
\rput(0,-2.8){{$\vdots$}} \rput(1.6,-2.8){{$\vdots$}} \rput(2.8,-2.8){{$\ddots$}}
\rput(4.2,-2.8){{$\vdots$}} \rput(6,-2.8){{$\vdots$}} \rput(4.9,-.8){{or}}
\end{pspicture}
$$
\caption{A proper total coloring distinguishing adjacent vertices
by sums of the product graph $P_m \times P_n$ where $m \ge 3$.} \label{pnpmfig}
\end{figure}

Last, if $m, n\ge 3$, then the maximum degree $\Delta$ of $P_m \times P_n$ is $4$.
If there exists only one vertex of the maximum degree,
then it must be $P_3 \times P_3$ which was handled previously.
Now, we assume there exist two adjacent vertices of the maximum degree.
The method we are going to use is sorely depends on the parity of $i$.
We define a coloring $\psi$ on $P_m \times P_n$ by

$$ \psi ((P_m \times P_n)_i) = \left\{ \begin{array}{cl}
\phi((P_m \times P_2)_1)& \rm{if} ~\textit{i}~\rm{is}~\rm{odd}, \\
\phi((P_m \times P_2)_2) & \rm{if} ~\textit{i}~\rm{is}~\rm{even}, \end{array}\right.$$
$$ \psi (h((P_m \times P_n)_i)) = \left\{ \begin{array}{cl}
5 & {\rm{if}}~i~{\rm{is}}~{\rm{odd}}, \\
6 & {\rm{if}}~i~{\rm{is}}~{\rm{even}}, \end{array}\right.$$
as illustrated in Figure~\ref{pnpmfig}.
One can easily check that $\psi$ is a proper total $6$ coloring
distinguishing adjacent vertices by sums of the graph $P_m \times P_n$,
which means $\tndi (P_m \times P_n) \le \Delta +2$.
By Theorem~\ref{obs1}, the opposite inequality $\Delta +2 \le\tndi (P_m \times P_n)$ holds.
Therefore, it completes the proof of theorem.
\end{proof}
\end{theorem}

\begin{figure}
$$
\begin{pspicture}[shift=-2.2](-3,-3)(3,3.2)
\psline(0,3)(-3,-1.8)(3,-1.8)(0,3)
\psline(0,1.5)(-1.5,-1)(1.5,-1)(0,1.5)
\psline(0,3)(0,1.5)
\psline(-3,-1.8)(-1.5,-1)
\psline(3,-1.8)(1.5,-1)
\pscircle[fillcolor=lightgray, fillstyle=solid, linewidth=1pt](0,3){.3}
\pscircle[fillcolor=lightgray, fillstyle=solid, linewidth=1pt](-3,-1.8){.3}
\pscircle[fillcolor=lightgray, fillstyle=solid, linewidth=1pt](3,-1.8){0.3}
\pscircle[fillcolor=lightgray, fillstyle=solid, linewidth=1pt](0,1.5){.3}
\pscircle[fillcolor=lightgray, fillstyle=solid, linewidth=1pt](-1.5,-1){.3}
\pscircle[fillcolor=lightgray, fillstyle=solid, linewidth=1pt](1.5,-1){.3}
\rput(0,3){{$1$}}\rput(-3,-1.8){{$3$}}\rput(3,-1.8){{$2$}}
\rput(0,1.5){{$2$}}\rput(-1.5,-1){{$1$}}\rput(1.5,-1){{$3$}}
\rput(2.2;87){{\footnotesize{3}}}\rput(2.4;207){{\footnotesize{4}}}
\rput(2.4;332){{\footnotesize{1}}}\rput(1;20){{\footnotesize{4}}}
\rput(1.8;20){{\footnotesize{4}}}\rput(1;160){{\footnotesize{5}}}
\rput(1.8;160){{\footnotesize{2}}} \rput(2.8;270){{$(a)$}}
\rput(1.2;270){{\footnotesize{2}}}\rput(2;270){{\footnotesize{5}}}
\end{pspicture}
\quad\quad
\begin{pspicture}[shift=-2](-3,-3)(3,3.2)
\psline(2.75;-67.5)(3;-45)(3;0)(3;45)(3;90)(3;135)(3;180)(3;225)(2.75;247.5)
\psline(1.35;-67.5)(1.5;-45)(1.5;0)(1.5;45)(1.5;90)(1.5;135)(1.5;180)(1.5;225)(1.35;247.5)
\psline(3;-45)(1.5;-45) \psline(3;0)(1.5;0) \psline(3;45)(1.5;45) \psline(3;225)(1.5;225)
\psline(3;90)(1.5;90) \psline(3;135)(1.5;135) \psline(3;180)(1.5;180)
\pscircle[fillcolor=lightgray, fillstyle=solid, linewidth=1pt](1.5;-45){.3}
\pscircle[fillcolor=lightgray, fillstyle=solid, linewidth=1pt](1.5;0){.3}
\pscircle[fillcolor=lightgray, fillstyle=solid, linewidth=1pt](1.5;45){0.3}
\pscircle[fillcolor=lightgray, fillstyle=solid, linewidth=1pt](1.5;90){.3}
\pscircle[fillcolor=lightgray, fillstyle=solid, linewidth=1pt](1.5;135){.3}
\pscircle[fillcolor=lightgray, fillstyle=solid, linewidth=1pt](1.5;180){.3}
\pscircle[fillcolor=lightgray, fillstyle=solid, linewidth=1pt](1.5;225){.3}
\pscircle[fillcolor=lightgray, fillstyle=solid, linewidth=1pt](3;-45){.3}
\pscircle[fillcolor=lightgray, fillstyle=solid, linewidth=1pt](3;0){.3}
\pscircle[fillcolor=lightgray, fillstyle=solid, linewidth=1pt](3;45){0.3}
\pscircle[fillcolor=lightgray, fillstyle=solid, linewidth=1pt](3;90){.3}
\pscircle[fillcolor=lightgray, fillstyle=solid, linewidth=1pt](3;135){.3}
\pscircle[fillcolor=lightgray, fillstyle=solid, linewidth=1pt](3;180){.3}
\pscircle[fillcolor=lightgray, fillstyle=solid, linewidth=1pt](3;225){.3}
\rput(3;-45){{$3$}}\rput(3;0){{$4$}}\rput(3;45){{$1$}}
\rput(3;90){{$3$}}\rput(3;135){{$1$}}\rput(3;180){{$3$}}\rput(3;225){{$1$}}
\rput(1.5;-45){{$1$}}\rput(1.5;0){{$2$}}\rput(1.5;45){{$3$}}
\rput(1.5;90){{$1$}}\rput(1.5;135){{$3$}}\rput(1.5;180){{$1$}}\rput(1.5;225){{$3$}}
\rput(2.95;-67.5){{\footnotesize{2}}}\rput(2.95;-22.5){{\footnotesize{1}}}
\rput(2.95;22.5){{\footnotesize{3}}}\rput(2.95;67.5){{\footnotesize{2}}}
\rput(2.95;112.5){{\footnotesize{4}}}\rput(2.95;157.5){{\footnotesize{2}}}
\rput(2.95;202.5){{\footnotesize{4}}}\rput(2.95;247.5){{\footnotesize{2}}}
\rput(1.2;-67.5){{\footnotesize{4}}}\rput(1.2;-22.5){{\footnotesize{3}}}
\rput(1.2;22.5){{\footnotesize{1}}}\rput(1.2;67.5){{\footnotesize{4}}}
\rput(1.2;112.5){{\footnotesize{2}}}\rput(1.2;157.5){{\footnotesize{4}}}
\rput(1.2;202.5){{\footnotesize{2}}}\rput(1.2;247.5){{\footnotesize{4}}}
\rput(2.25;-49){{\footnotesize{5}}}\rput(2.25;-4){{\footnotesize{5}}}
\rput(2.25;41){{\footnotesize{5}}}\rput(2.25;86){{\footnotesize{5}}}
\rput(2.25;131){{\footnotesize{5}}}\rput(2.25;176){{\footnotesize{5}}}
\rput(2.25;221){{\footnotesize{5}}} \rput(2.25;270){{$\cdots$}}  \rput(3;270){{$(b)$}}
\end{pspicture}
$$
\caption{A proper total coloring distinguishing adjacent vertices
by sums of the product graph $(a)$ $C_3 \times P_2$ and $(b)$ $C_{2k+1} \times P_2$, $k \ge 2$. } \label{c3p2fig}
\end{figure}

\begin{theorem}
Let $C_m$ be a cycle of order $m \ge3$ and $P_n$ be a path of order $n \ge 2$ .
Then $\tndi (C_m \times P_n) = \Delta +2$.
\label{pncm}
\begin{proof}
These $C_m \times P_n$ graphs always have two adjacent vertices $x$, $y$ such that $deg(x)=deg(y)=\Delta$. By
Theorem~\ref{obs1}, we obtain the inequality $\tndi (C_m \times P_n) \ge \Delta +2$.

First look at the case $n=2$. Since $\Delta =3$, to show $C_m \times P_2$ is $\tndi$ class II,
we need to find a proper total $5$ coloring distinguishing adjacent vertices by sums.
If $m$ is even, then one can see that $C_m \times P_2$ is regular bipartite.
By Theorem~\ref{tndioldthm} (6), it is $\tndi$ class II.
To handle the remaining general cases, let us find a coloring of $C_m \times P_2$ as follows.

If $m$ is odd, let us denote it by $m=2k+1$. For $m=3$, we color it as depicted in Figure~\ref{c3p2fig} $(a)$.
If $2k+1 \ge 5$, we color the vertices of
$(C_{2k+1} \times P_2)_1$ with colors $4, 1, 3, 1, 3, \ldots, 1, 3$
and the edges with $1, 3, 2, 4, 2, 4, \ldots, 4, 2$.
Similarly, color the consecutive vertices of $(C_{2k+1} \times P_2)_2$ with $2, 3, 1, 3, 1, \ldots ,3 ,1 $
and consecutive edges with $3, 1, 4, 2, 4, 2, \ldots, 4$.
Color the reminding edges with the color $5$ as illustrated in Figure~\ref{c3p2fig} $(b)$.

If $m=2k\ge 4$, we color the vertices of
$(C_{2k} \times P_2)_1$ with colors $ 1, 3, 1, 3, \ldots, 1, 3$
and the edges with $ 2, 4, 2, 4, \ldots, 2, 4$.
Similarly, color the consecutive vertices of $(C_{2k}
\times P_2)_2$ with $ 3, 1, 3, 1, \ldots ,3 ,1 $
and consecutive edges with $4, 2, 4, 2, \ldots, 4,2$.
Color the reminder edges with the color $5$ then this is a proper
total $5$ coloring distinguishing adjacent vertices by sums of $C_{2k} \times P_2$.

Now we are going to deal with the case $n \ge 3$.
Since the maximum degree $\Delta$ of $C_m \times P_n$ is $4$,
to show $C_m \times P_n$ is $\tndi$ class II,
we need to find a proper
total $6$ coloring distinguishing adjacent vertices by sums.
The above coloring for $C_{m} \times P_2$ is denoted by $\phi$.

\begin{figure}
$$
\begin{pspicture}[shift=-3](-2.7,-3)(2.7,1.8)
\psline[linecolor=darkred, linewidth=1.5pt](-2.7,1.6)(-2.7,-2.4)
\psline[linecolor=darkred, linestyle=dotted, linewidth=1.5pt](-2.7,-2.4)(-2.7,-3)
\psline[linecolor=darkred, arrowscale=4]{->}(-2.7,0.56)(-2.7,0.72)
\psline[linecolor=darkred, linewidth=1.5pt](2.7,1.6)(2.7,-2.4)
\psline[linecolor=darkred, linestyle=dotted, linewidth=1.5pt](2.7,-2.4)(2.7,-3)
\psline[linecolor=darkred, arrowscale=4]{->}(2.7,0.56)(2.7,0.72)
\psline(-2.7,1.6)(2,1.6)
\psline(1.6,1.6)(1.6,-2.4)
\psline(-2.7,-1.6)(2,-1.6)
\psline(-2.7,0)(2,0)
\psline(-1.6,1.6)(-1.6,-2.4)
\psline(0,1.6)(0,-2.4)
\psline(2,0)(2.7,0)
\psline(2,-1.6)(2.7,-1.6)
\psline(2,1.6)(2.7,1.6)
\pscircle[fillcolor=lightgray, fillstyle=solid, linewidth=1pt](-1.6,1.6){.3}
\pscircle[fillcolor=lightgray, fillstyle=solid, linewidth=1pt](0,1.6){.3}
\pscircle[fillcolor=lightgray, fillstyle=solid, linewidth=1pt](1.6,1.6){.3}
\pscircle[fillcolor=lightgray, fillstyle=solid, linewidth=1pt](-1.6,0){.3}
\pscircle[fillcolor=lightgray, fillstyle=solid, linewidth=1pt](0,0){.3}
\pscircle[fillcolor=lightgray, fillstyle=solid, linewidth=1pt](1.6,0){.3}
\pscircle[fillcolor=lightgray, fillstyle=solid, linewidth=1pt](-1.6,-1.6){.3}
\pscircle[fillcolor=lightgray, fillstyle=solid, linewidth=1pt](0,-1.6){.3}
\pscircle[fillcolor=lightgray, fillstyle=solid, linewidth=1pt](1.6,-1.6){.3}
\rput(-1.6,1.6){{$1$}} \rput(0,1.6){{$3$}} \rput(1.6,1.6){{$2$}}
\rput(-1.6,0){{$3$}} \rput(0,0){{$2$}} \rput(1.6,0){{$1$}}
\rput(-1.6,-1.6){{$1$}}\rput(0,-1.6){{$3$}}\rput(1.6,-1.6){{$2$}}
\rput(1.76,0.8){{\footnotesize{3}}}\rput(0.16,0.8){{\footnotesize{1}}}
\rput(-1.44,0.8){{\footnotesize{4}}}
\rput(1.76,-0.8){{\footnotesize{6}}}\rput(0.16,-0.8){{\footnotesize{6}}}
\rput(-1.44,-0.8){{\footnotesize{6}}}
\rput(1.76,-2.3){{\footnotesize{3}}}\rput(0.16,-2.3){{\footnotesize{1}}}
\rput(-1.44,-2.3){{\footnotesize{4}}}
\rput(0.8,1.85){{\footnotesize{4}}}\rput(0.8,0.25){{\footnotesize{4}}}
\rput(0.8,-1.35){{\footnotesize{1}}}
\rput(-0.8,1.85){{\footnotesize{2}}}\rput(-0.8,0.25){{\footnotesize{5}}}
\rput(-0.8,-1.35){{\footnotesize{2}}}
\rput(-2.15,1.85){{\footnotesize{5}}}\rput(-2.15,0.25){{\footnotesize{2}}}
\rput(-2.15,-1.35){{\footnotesize{5}}}
\rput(2.15,1.85){{\footnotesize{5}}}\rput(2.15,0.25){{\footnotesize{2}}}
\rput(2.15,-1.35){{\footnotesize{5}}}
\rput(-1.6,-2.7){{$\vdots$}}\rput(0,-2.7){{$\vdots$}}
\rput(1.6,-2.7){{$\vdots$}}
\end{pspicture}
$$
\caption{A proper total coloring distinguishing adjacent vertices
by sums of the product graph $C_3 \times P_n$.} \label{c3pnfig}
\end{figure}
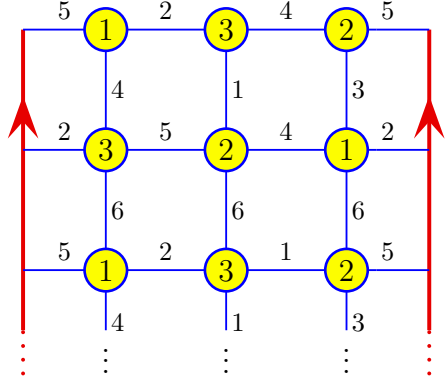

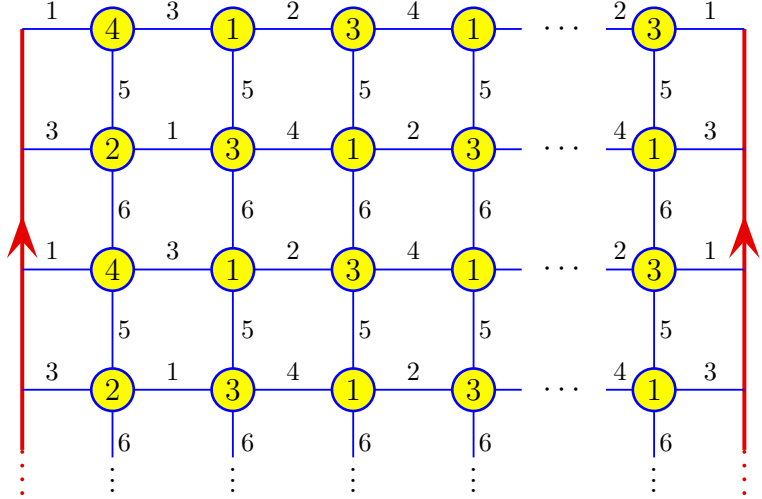
\begin{figure}
$$
\begin{pspicture}[shift=-3](-2.8,-3)(6.8,3.5)
\psline[linecolor=darkred, linewidth=1.5pt](-2.8,3.2)(-2.8,-2.4)
\psline[linecolor=darkred, linestyle=dotted, linewidth=1.5pt](-2.8,-2.4)(-2.8,-3)
\psline[linecolor=darkred, arrowscale=4]{->}(-2.8,0.56)(-2.8,0.72)
\psline[linecolor=darkred, linewidth=1.5pt](6.8,3.2)(6.8,-2.4)
\psline[linecolor=darkred, linestyle=dotted, linewidth=1.5pt](6.8,-2.4)(6.8,-3)
\psline[linecolor=darkred, arrowscale=4]{->}(6.8,0.56)(6.8,0.72)
\psline(-2.8,-1.6)(3.84,-1.6) \psline(4.96,-1.6)(6.8,-1.6)
\psline(-2.8,0)(3.84,0) \psline(4.96,0)(6.8,0)
\psline(-2.8,1.6)(3.84,1.6) \psline(4.96,1.6)(6.8,1.6)
\psline(-2.8,3.2)(3.84,3.2) \psline(4.96,3.2)(6.8,3.2)
\psline(0,3.2)(0,-2.5) \psline(1.6,3.2)(1.6,-2.5)
\psline(-1.6,3.2)(-1.6,-2.5) \psline(3.2,3.2)(3.2,-2.5) \psline(5.6,3.2)(5.6,-2.5)
\pscircle[fillcolor=lightgray, fillstyle=solid, linewidth=1pt](-1.6,3.2){.3}
\pscircle[fillcolor=lightgray, fillstyle=solid, linewidth=1pt](0,3.2){.3}
\pscircle[fillcolor=lightgray, fillstyle=solid, linewidth=1pt](1.6,3.2){.3}
\pscircle[fillcolor=lightgray, fillstyle=solid, linewidth=1pt](3.2,3.2){.3}
\pscircle[fillcolor=lightgray, fillstyle=solid, linewidth=1pt](5.6,3.2){.3}
\pscircle[fillcolor=lightgray, fillstyle=solid, linewidth=1pt](-1.6,1.6){.3}
\pscircle[fillcolor=lightgray, fillstyle=solid, linewidth=1pt](0,1.6){.3}
\pscircle[fillcolor=lightgray, fillstyle=solid, linewidth=1pt](1.6,1.6){.3}
\pscircle[fillcolor=lightgray, fillstyle=solid, linewidth=1pt](3.2,1.6){.3}
\pscircle[fillcolor=lightgray, fillstyle=solid, linewidth=1pt](5.6,1.6){.3}
\pscircle[fillcolor=lightgray, fillstyle=solid, linewidth=1pt](3.2,0){.3}
\pscircle[fillcolor=lightgray, fillstyle=solid, linewidth=1pt](-1.6,0){.3}
\pscircle[fillcolor=lightgray, fillstyle=solid, linewidth=1pt](0,0){.3}
\pscircle[fillcolor=lightgray, fillstyle=solid, linewidth=1pt](1.6,0){.3}
\pscircle[fillcolor=lightgray, fillstyle=solid, linewidth=1pt](5.6,0){.3}
\pscircle[fillcolor=lightgray, fillstyle=solid, linewidth=1pt](3.2,-1.6){.3}
\pscircle[fillcolor=lightgray, fillstyle=solid, linewidth=1pt](-1.6,-1.6){.3}
\pscircle[fillcolor=lightgray, fillstyle=solid, linewidth=1pt](0,-1.6){.3}
\pscircle[fillcolor=lightgray, fillstyle=solid, linewidth=1pt](1.6,-1.6){.3}
\pscircle[fillcolor=lightgray, fillstyle=solid, linewidth=1pt](5.6,-1.6){.3}
\rput(-1.6,3.2){{$4$}}\rput(0,3.2){{$1$}}\rput(1.6,3.2){{$3$}}\rput(3.2,3.2){{$1$}}
\rput(-1.6,1.6){{$2$}}\rput(0,1.6){{$3$}}\rput(1.6,1.6){{$1$}}\rput(3.2,1.6){{$3$}}
\rput(-1.6,0){{$4$}}\rput(0,0){{$1$}}\rput(1.6,0){{$3$}}\rput(3.2,0){{$1$}}
\rput(-1.6,-1.6){{$2$}}\rput(0,-1.6){{$3$}}\rput(1.6,-1.6){{$1$}}\rput(3.2,-1.6){{$3$}}
\rput(5.6,3.2){{$3$}}\rput(5.6,1.6){{$1$}}\rput(5.6,0){{$3$}} \rput(5.6,-1.6){{$1$}}
\rput(-2.4,3.45){{\footnotesize{1}}}\rput(-0.8,3.45){{\footnotesize{3}}}
\rput(0.8,3.45){{\footnotesize{2}}}\rput(2.4,3.45){{\footnotesize{4}}}
\rput(-2.4,1.85){{\footnotesize{3}}}\rput(-0.8,1.85){{\footnotesize{1}}}
\rput(0.8,1.85){{\footnotesize{4}}}\rput(2.4,1.85){{\footnotesize{2}}}
\rput(-2.4,0.25){{\footnotesize{1}}}\rput(-0.8,0.25){{\footnotesize{3}}}
\rput(0.8,0.25){{\footnotesize{2}}}\rput(2.4,0.25){{\footnotesize{4}}}
\rput(-2.4,-1.35){{\footnotesize{3}}}\rput(-0.8,-1.35){{\footnotesize{1}}}
\rput(0.8,-1.35){{\footnotesize{4}}}\rput(2.4,-1.35){{\footnotesize{2}}}
\rput(-1.44,2.4){{\footnotesize{5}}}\rput(-1.44,0.8){{\footnotesize{6}}}
\rput(-1.44,-0.8){{\footnotesize{5}}}\rput(0.2,2.4){{\footnotesize{5}}}
\rput(0.2,0.8){{\footnotesize{6}}}\rput(0.2,-0.8){{\footnotesize{5}}}
\rput(1.76,2.4){{\footnotesize{5}}}\rput(1.76,0.8){{\footnotesize{6}}}
\rput(1.76,-0.8){{\footnotesize{5}}}\rput(3.36,2.4){{\footnotesize{5}}}
\rput(3.36,0.8){{\footnotesize{6}}}\rput(3.36,-0.8){{\footnotesize{5}}}
\rput(5.76,2.4){{\footnotesize{5}}}\rput(5.76,0.8){{\footnotesize{6}}}
\rput(5.76,-0.8){{\footnotesize{5}}}\rput(-1.44,-2.3){{\footnotesize{6}}}
\rput(5.15,3.45){{\footnotesize{2}}}\rput(5.15,1.85){{\footnotesize{4}}}
\rput(5.15,0.25){{\footnotesize{2}}}\rput(5.15,-1.35){{\footnotesize{4}}}
\rput(6.35,3.45){{\footnotesize{1}}}\rput(6.35,1.85){{\footnotesize{3}}}
\rput(6.35,0.25){{\footnotesize{1}}}\rput(6.35,-1.35){{\footnotesize{3}}}
\rput(1.76,-2.3){{\footnotesize{6}}}\rput(0.2,-2.3){{\footnotesize{6}}}
\rput(5.76,-2.3){{\footnotesize{6}}}\rput(3.36,-2.3){{\footnotesize{6}}}
\rput(-1.6,-2.7){{$\vdots$}}\rput(0,-2.7){{$\vdots$}}
\rput(1.6,-2.7){{$\vdots$}} \rput(3.2,-2.7){{$\vdots$}}\rput(5.6,-2.7){{$\vdots$}}
\rput(4.4,3.2){{$\cdots$}} \rput(4.4,1.6){{$\cdots$}}
\rput(4.4,0){{$\cdots$}} \rput(4.4,-1.6){{$\cdots$}}
\end{pspicture}
$$
\caption{A proper total coloring distinguishing adjacent vertices
by sums of the product graph $C_m \times P_n$ where $m$ is bigger than $3$ and odd.} \label{cmpnfig1}
\end{figure}

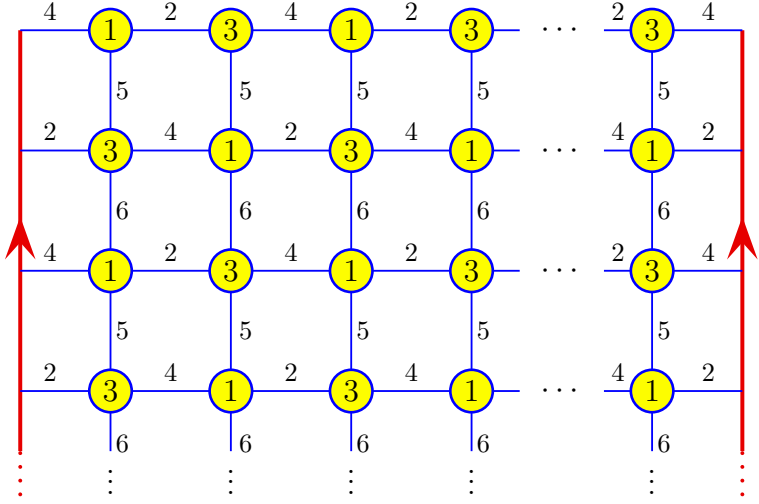
\begin{figure}
$$
\begin{pspicture}[shift=-3](-2.8,-3)(6.8,3.5)
\psline[linecolor=darkred, linewidth=1.5pt](-2.8,3.2)(-2.8,-2.4)
\psline[linecolor=darkred, linestyle=dotted, linewidth=1.5pt](-2.8,-2.4)(-2.8,-3)
\psline[linecolor=darkred, arrowscale=4]{->}(-2.8,0.56)(-2.8,0.72)
\psline[linecolor=darkred, linewidth=1.5pt](6.8,3.2)(6.8,-2.4)
\psline[linecolor=darkred, linestyle=dotted, linewidth=1.5pt](6.8,-2.4)(6.8,-3)
\psline[linecolor=darkred, arrowscale=4]{->}(6.8,0.56)(6.8,0.72)
\psline(-2.8,-1.6)(3.84,-1.6) \psline(4.96,-1.6)(6.8,-1.6)
\psline(-2.8,0)(3.84,0) \psline(4.96,0)(6.8,0)
\psline(-2.8,1.6)(3.84,1.6) \psline(4.96,1.6)(6.8,1.6)
\psline(-2.8,3.2)(3.84,3.2) \psline(4.96,3.2)(6.8,3.2)
\psline(0,3.2)(0,-2.4) \psline(1.6,3.2)(1.6,-2.4)
\psline(-1.6,3.2)(-1.6,-2.4) \psline(3.2,3.2)(3.2,-2.4) \psline(5.6,3.2)(5.6,-2.4)
\pscircle[fillcolor=lightgray, fillstyle=solid, linewidth=1pt](-1.6,3.2){.3}
\pscircle[fillcolor=lightgray, fillstyle=solid, linewidth=1pt](0,3.2){.3}
\pscircle[fillcolor=lightgray, fillstyle=solid, linewidth=1pt](1.6,3.2){.3}
\pscircle[fillcolor=lightgray, fillstyle=solid, linewidth=1pt](3.2,3.2){.3}
\pscircle[fillcolor=lightgray, fillstyle=solid, linewidth=1pt](5.6,3.2){.3}
\pscircle[fillcolor=lightgray, fillstyle=solid, linewidth=1pt](-1.6,1.6){.3}
\pscircle[fillcolor=lightgray, fillstyle=solid, linewidth=1pt](0,1.6){.3}
\pscircle[fillcolor=lightgray, fillstyle=solid, linewidth=1pt](1.6,1.6){.3}
\pscircle[fillcolor=lightgray, fillstyle=solid, linewidth=1pt](3.2,1.6){.3}
\pscircle[fillcolor=lightgray, fillstyle=solid, linewidth=1pt](5.6,1.6){.3}
\pscircle[fillcolor=lightgray, fillstyle=solid, linewidth=1pt](3.2,0){.3}
\pscircle[fillcolor=lightgray, fillstyle=solid, linewidth=1pt](-1.6,0){.3}
\pscircle[fillcolor=lightgray, fillstyle=solid, linewidth=1pt](0,0){.3}
\pscircle[fillcolor=lightgray, fillstyle=solid, linewidth=1pt](1.6,0){.3}
\pscircle[fillcolor=lightgray, fillstyle=solid, linewidth=1pt](5.6,0){.3}
\pscircle[fillcolor=lightgray, fillstyle=solid, linewidth=1pt](3.2,-1.6){.3}
\pscircle[fillcolor=lightgray, fillstyle=solid, linewidth=1pt](-1.6,-1.6){.3}
\pscircle[fillcolor=lightgray, fillstyle=solid, linewidth=1pt](0,-1.6){.3}
\pscircle[fillcolor=lightgray, fillstyle=solid, linewidth=1pt](1.6,-1.6){.3}
\pscircle[fillcolor=lightgray, fillstyle=solid, linewidth=1pt](5.6,-1.6){.3}
\rput(-1.6,3.2){{$1$}}\rput(0,3.2){{$3$}}\rput(1.6,3.2){{$1$}}\rput(3.2,3.2){{$3$}}
\rput(-1.6,1.6){{$3$}}\rput(0,1.6){{$1$}}\rput(1.6,1.6){{$3$}}\rput(3.2,1.6){{$1$}}
\rput(-1.6,0){{$1$}}\rput(0,0){{$3$}}\rput(1.6,0){{$1$}}\rput(3.2,0){{$3$}}
\rput(-1.6,-1.6){{$3$}}\rput(0,-1.6){{$1$}}\rput(1.6,-1.6){{$3$}}\rput(3.2,-1.6){{$1$}}
\rput(5.6,3.2){{$3$}}\rput(5.6,1.6){{$1$}}\rput(5.6,0){{$3$}}\rput(5.6,-1.6){{$1$}}
\rput(-2.4,3.45){{\footnotesize{4}}} \rput(-0.8,3.45){{\footnotesize{2}}}
\rput(0.8,3.45){{\footnotesize{4}}} \rput(2.4,3.45){{\footnotesize{2}}}
\rput(-2.4,1.85){{\footnotesize{2}}} \rput(-0.8,1.85){{\footnotesize{4}}}
\rput(0.8,1.85){{\footnotesize{2}}} \rput(2.4,1.85){{\footnotesize{4}}}
\rput(-2.4,0.25){{\footnotesize{4}}} \rput(-0.8,0.25){{\footnotesize{2}}}
\rput(0.8,0.25){{\footnotesize{4}}} \rput(2.4,0.25){{\footnotesize{2}}}
\rput(-2.4,-1.35){{\footnotesize{2}}} \rput(-0.8,-1.35){{\footnotesize{4}}}
\rput(0.8,-1.35){{\footnotesize{2}}}\rput(2.4,-1.35){{\footnotesize{4}}}
\rput(-1.44,-0.8){{\footnotesize{5}}} \rput(0.2,-0.8){{\footnotesize{5}}}
\rput(1.76,-0.8){{\footnotesize{5}}} \rput(3.36,-0.8){{\footnotesize{5}}}
\rput(5.76,-0.8){{\footnotesize{5}}} \rput(-1.44,0.8){{\footnotesize{6}}}
\rput(0.2,0.8){{\footnotesize{6}}} \rput(1.76,0.8){{\footnotesize{6}}}
\rput(3.36,0.8){{\footnotesize{6}}} \rput(5.76,0.8){{\footnotesize{6}}}
\rput(-1.44,2.4){{\footnotesize{5}}} \rput(0.2,2.4){{\footnotesize{5}}}
\rput(1.76,2.4){{\footnotesize{5}}} \rput(3.36,2.4){{\footnotesize{5}}}
\rput(5.76,2.4){{\footnotesize{5}}} \rput(-1.44,-2.3){{\footnotesize{6}}}
\rput(0.2,-2.3){{\footnotesize{6}}} \rput(1.76,-2.3){{\footnotesize{6}}}
\rput(6.35,3.45){{\footnotesize{4}}}\rput(6.35,1.85){{\footnotesize{2}}}
\rput(6.35,0.25){{\footnotesize{4}}}\rput(6.35,-1.35){{\footnotesize{2}}}
\rput(3.36,-2.3){{\footnotesize{6}}} \rput(5.76,-2.3){{\footnotesize{6}}}
\rput(5.15,3.45){{\footnotesize{2}}} \rput(5.15,1.85){{\footnotesize{4}}}
\rput(5.15,0.25){{\footnotesize{2}}} \rput(5.15,-1.35){{\footnotesize{4}}}
\rput(-1.6,-2.7){{$\vdots$}}\rput(0,-2.7){{$\vdots$}}
\rput(1.6,-2.7){{$\vdots$}} \rput(3.2,-2.7){{$\vdots$}}
\rput(5.6,-2.7){{$\vdots$}}
\rput(4.4,3.2){{$\cdots$}} \rput(4.4,1.6){{$\cdots$}}
\rput(4.4,0){{$\cdots$}} \rput(4.4,-1.6){{$\cdots$}}
\end{pspicture}
$$
\caption{A proper total coloring distinguishing adjacent vertices by sums
of the product graph $C_m \times P_n$ where $m$ is bigger than $3$ and even.} \label{cmpnfig2}
\end{figure}

If $m=3$, then color the consecutive vertices of $(C_3 \times P_n)_1$ with $1, 3, 2$ and
corresponding consecutive edges with $2, 4, 5$. Similarly, we color the consecutive
vertices of $(C_3 \times P_n)_2$ with $3, 2, 1$ and corresponding consecutive edges with $5, 4, 2$ as 
we have used for $C_3 \times P_2$ as illustrated in Figure~\ref{c3p2fig} $(a)$.
We expand these colorings on $(C_3 \times P_n)_1$ and $(C_3 \times P_n)_2$
to a coloring $\psi$ of $(C_3 \times P_n)$ by

$$ \psi ((C_3 \times P_n)_i) = \left\{ \begin{array}{cl}
\phi((C_3 \times P_2)_1)& {\rm{if}} ~i~{\rm{is}}~{\rm{odd}}, \\
\phi((C_3 \times P_2)_2)& {\rm{if}} ~i~{\rm{is}}~{\rm{even}}, \end{array}\right.$$
$$ \psi (h((C_3 \times P_n)_i)) = \left\{ \begin{array}{cl}
4, 1, 3
& {\rm{if}} ~i~{\rm{is}}~{\rm{odd}}, \\  6, 6, 6 & \rm{if}
~\textit{i}~\rm{is}~\rm{even}, \end{array}\right.$$
as illustrated in Figure~\ref{c3pnfig}.

If $m=2k+1 \ge 5$, we define a coloring $\psi$ by

$$ \psi ((C_{2k+1} \times P_n)_i) = \left\{ \begin{array}{cl}
\phi((C_{2k+1} \times P_2)_1)& \rm{if} ~\textit{i}~\rm{is}~\rm{odd}, \\
\phi((C_{2k+1} \times P_2)_2) & \rm{if} ~\textit{i}~\rm{is}~\rm{even}, \end{array}\right.$$
$$ \psi (h((C_{2k+1} \times P_n)_i)) = \left\{ \begin{array}{cl}
5 & {\rm{if}} ~i~{\rm{is}}~{\rm{odd}}, \\
6 & {\rm{if}} ~i~{\rm{is}}~{\rm{even}}, \end{array}\right.$$
as depicted in Figure~\ref{cmpnfig1}.

If $m=2k \ge 4$, we define a coloring $\psi$ by

$$ \psi ((C_{2k} \times P_n)_i) =
\left\{ \begin{array}{cl}
\phi((C_{2k} \times P_2)_1)& \rm{if} ~\textit{i}~\rm{is}~\rm{odd}, \\
\phi((C_{2k} \times P_2)_2) & \rm{if} ~\textit{i}~\rm{is}~\rm{even}, \end{array}\right.$$
$$ \psi (h((C_{2k} \times P_n)_i)) = \left\{ \begin{array}{cl}
5 & {\rm{if}} ~i~{\rm{is}}~{\rm{odd}}, \\
6 & {\rm{if}} ~i~{\rm{is}}~{\rm{even}}, \end{array}\right.$$
as illustrated in Figure~\ref{cmpnfig2}.

It is not difficult check that these three colorings $\psi$'s are proper total
$6$ colorings distinguishing adjacent vertices by sums.
Therefore, it completes the proof of theorem.
\end{proof}
\end{theorem}

\begin{theorem}
Let $S_m$ be the star graph of order $m \ge 2$ and $P_n$ be the path of order $n\ge 2$. Then
$$\tndi (S_m \times P_n)= \left\{ \begin{array}{cl}
\Delta +1 & {\rm{if}} ~n=3, \\
\Delta +2 & {\rm{if}} ~n=2~{\rm{or}}~n \ge 4. \end{array}\right.$$
\label{pnsm}
\begin{proof}
First let us deal with the case $n=2$.
Then it has two adjacent vertices of the maximum degree $\Delta=m+1$.
To prove $S_m \times P_n$ is $\tndi$ class II, we need to find a proper total $m+3(=\Delta + 2)$
coloring distinguishing adjacent vertices by sums on $S_m \times P_2$
by Theorem~\ref{obs1}.

On the other hand, there exists a proper total $m+1$ coloring distinguishing
adjacent vertices by sums for $S_m$ by Theorem~\ref{tndioldthm} (3) denoted by $c$.
Define a new coloring $c'$ by $c'= c+1$.
Color $(S_m \times P_2)_1$ by $c$ and $(S_m \times P_2)_2$ by $c'$
and $h((S_m \times P_2)_1)=m+3$. Then one can easily check that this is a proper total
$(m+3)$ coloring distinguishing adjacent vertices by sums of $S_m \times P_2$.
This coloring for $S_{m} \times P_2$ is denoted by $\phi$.

Second, we assume $n \ge 4$, then $\Delta+2=m+4$ and $S_m \times P_n$ has two adjacent vertices of
the maximum degree. By Theorem~\ref{obs1}, we need to find
a proper total $(m+4)$ coloring distinguishing adjacent vertices by sums.
We define a coloring $\psi$
on $(S_m \times P_n)$ as follows,
$$ \psi ((S_m \times P_n)_i) = \left\{ \begin{array}{cl}
\phi((S_m \times P_2)_1)& \rm{if} ~\textit{i}~\rm{is}~\rm{odd}, \\
\phi((S_m \times P_2)_2) & \rm{if} ~\textit{i}~\rm{is}~\rm{even}, \end{array}\right.$$
$$ \psi (h((S_m \times P_n)_i)) = \left\{ \begin{array}{cl}
m+3 & {\rm{if}} ~i~{\rm{is}}~{\rm{odd}}, \\
m+4 & {\rm{if}} ~i~{\rm{is}}~{\rm{even}}. \end{array}\right.$$
One may easily check that $\psi$ is a proper total $m+4$ coloring distinguishing adjacent
vertices by sums. Therefore, $(S_m \times P_n)$ is $\tndi$ class II.

We are discussing he last case $n=3$. For this case, we get just one vertex of the maximum degree.
As we have seen for $P_3 \times P_3$ in Figure~\ref{p2pmfig} $(b)$,
$S_2 \times P_3=P_3 \times P_3$ is $\tndi$ class I.
For $m\ge 3$, we can show that $S_m \times P_3$ is $\tndi$ class I as follows.
Let $c$ be a proper total $m+1$ coloring distinguishing
adjacent vertices by sums of $S_m$,
and $c'\equiv c+1$ (mod $m+1$) be the second coloring provided in Theorem~\ref{snthm}.
We define a coloring $\psi$ on $S_m \times P_3$ as follows,
$$ \psi ((S_m \times P_3)_i) = \left\{ \begin{array}{cl}
c'((S_m \times P_2)_1)& {\rm{if}} ~i=1, 3, \\
c((S_m \times P_2)_2) & {\rm{if}} ~i=2, \end{array}\right.$$
$$ \psi (h((S_m \times P_3)_i)) = \left\{ \begin{array}{cl}
m+2 & {\rm{if}} ~i=1, \\
m+3 & {\rm{if}} ~i=2. \end{array}\right.$$
Then, one can check that $\psi$ is a proper total $m+3$ coloring distinguishing
adjacent vertices by sums for $S_m\times P_3$.
 It completes the proof of theorem.
\end{proof}
\end{theorem}

\begin{theorem}
Let $K_m$ be a
complete graph of order $m \ge 2$ and $P_n$ be a path of order $n\ge 2$.
Then $\tndi (K_m \times P_n) = \Delta+2$.
\label{pnkm}
\begin{proof}
Since $K_m \times P_n$ has two adjacent vertices of the maximum degree,
by Theorem~\ref{obs1}, to prove $K_m \times P_n$ is $\tndi$ class II,
we need find a proper total $\Delta +2$ coloring
distinguishing adjacent vertices by sums of $K_m \times P_n$.

First, let us assume that $n=2$ and $m$ is odd.
Then there exists a proper total coloring $c_0$
distinguishing adjacent vertices by sums of $K_m$ by Theorem~\ref{tndioldthm} (4).
This coloring $c_0$ originally found in~\cite{PW} has a very nice formula that for each vertex $v_i \in K_m$,
$$f(v_i)=  \frac{m^2+5m+6}{2} -(2i+1)$$
by the fact that $C(v_i)=\{1,2, \ldots, m+2\}
\setminus\{i,i+1\}$
as illustrated in Figure~\ref{k6fig} $(a)$ where
the vertex $v_i$ is colored $2i+1$ (mod $7$).
We define a coloring $\phi$ on $K_m \times P_2$ by coloring the $(K_m \times P_2)_2$ with the coloring $c_0$,
coloring the $(K_m \times P_2)_1$ with $c_0 '\equiv c_0 +1$ (mod $m+2$)
and $\phi(\{(v_i, w_1), (v_i, w_2)\})= g(i)+1$ for $i=1, 2, \ldots, m$.

If $n=2$ and $m$ is even, the coloring found in~\cite{PW} has a similar property that
$$f(v_i)=  \frac{m^2+3m+2}{2} -i ~\mod (m+1)$$
by the fact that $C(v_i)=\{1,2, \ldots, m+1\}
\setminus\{i\}$ as illustrated in Figure~\ref{k6fig} $(b)$
where the vertex $v_i$ is colored $2i$ (mod $7$).
We define a coloring $\phi$ on $K_m \times P_2$ by coloring $(K_m \times P_2)_2$ with the coloring $c_0$,
and by coloring $(K_m \times P_2)_1$ with $c_0 '\equiv c_0 +1$ (mod $m+1$)
and $\phi(\{(v_i, w_1), (v_i, w_2)\})= m+2$ for $i=1, 2, \ldots, m$.
Then one may easily check that this is a proper total coloring distinguishing adjacent vertices by sums.
The above coloring for $K_{m} \times P_2$ is denoted by $\phi$.

Let $m\ge 3$, then $\Delta+2=m+3$. If $m$ is odd, we define a coloring $\psi$ on $K_m \times P_n$ by
$$ \psi ((K_m \times P_n)_i) = \left\{ \begin{array}{cl}
\phi((K_m \times P_2)_1)& \rm{if} ~\textit{i}~\rm{is}~\rm{odd}, \\
\phi((K_m \times P_2)_2) & \rm{if} ~\textit{i}~\rm{is}~\rm{even}, \end{array}\right.$$
$$ \psi (h((K_m \times P_n)_i)) = \left\{ \begin{array}{cl}
m+3 & {\rm{if}} ~i~{\rm{is}}~{\rm{odd}}, \\
c_0 (h(K_m \times P_n)_i) & {\rm{if}} ~i~{\rm{is}}~{\rm{even}}. \end{array}\right.$$
 Clearly, $\psi$ is a proper total $m+3$ coloring distinguishing
 adjacent vertices by sums of $K_m \times P_n$.

 If $m$ is even, we define a coloring $\psi$ on $K_m \times P_n$ by
$$ \psi ((K_m \times P_n)_i) = \left\{ \begin{array}{cl}
\phi((K_m \times P_2)_1)& \rm{if} ~\textit{i}~\rm{is}~\rm{odd}, \\
\phi((K_m \times P_2)_2) & \rm{if} ~\textit{i}~\rm{is}~\rm{even}, \end{array}\right.$$
$$ \psi (h((K_m \times P_n)_i)) = \left\{ \begin{array}{cl}
m+2 & {\rm{if}} ~i~{\rm{is}}~{\rm{odd}}, \\
m+3 & {\rm{if}} ~i~{\rm{is}}~{\rm{even}}. \end{array}\right.$$
It is easy to see that $\psi$ is a proper total $m+3$ coloring distinguishing adjacent
vertices by sums of $K_m \times P_n$.
\end{proof}
\end{theorem}

\begin{theorem}
Let $W_m$ be a wheel graph of
order $m \ge 3$ and $P_n$ be a path of order $n\ge 2$. Then,
$$\tndi (W_m \times P_n)= \left\{ \begin{array}{cl}
\Delta +1 & {\rm{if}} ~n=3~{\rm{and}}~m\ge 4, \\
\Delta +2 & {\rm{Otherwise}}. \end{array}\right.$$
\label{pnwm}
\begin{proof}
Since $W_3=K_4$, $W_3 \times P_n$ is $\tndi$ class II as proved in Theorem~\ref{pnkm}.
If $n=2$ and $m \ge 4$, then the maximum degree is $m+1$.
There exists a proper total $m+1$ coloring $c$
distinguishing adjacent vertices by sums of $W_m$ by Theorem~\ref{wnthm}.
Let $c'=c+1$. Let us define a coloring $\phi$ of $W_m \times P_2$ by coloring
$(W_m \times P_2)_1$ with the coloring $c$, $(W_m \times P_2)_2$ with the coloring
$c'$ and $\phi(h(W_m \times P_2)_1)=m+3$. Then $\tndi (W_m \times P_2)=m+3$.

Next if $n\ge 4$ and $m\ge 4$, then maximum degree is $m+2$ and
there exist two adjacent vertices of the maximum degree.
We define a coloring $\psi$ on $W_m \times P_n$ as follows
$$ \psi((W_m \times P_n)_i) = \left\{ \begin{array}{cl}
\phi((W_m \times P_2)_1) & \rm{if} ~\textit{i}~\rm{is}~\rm{odd}, \\
\phi((W_m \times P_2)_2) & \rm{if} ~\textit{i}~\rm{is}~\rm{even}, \end{array}\right.$$
$$ \psi (h((W_m \times P_n)_i)) = \left\{ \begin{array}{cl}
m+3 & {\rm{if}} ~i~{\rm{is}~\rm{odd}}, \\
m+4 & {\rm{if}} ~i~{\rm{is}~\rm{even}}. \end{array}\right.$$
One may easily check that $\psi$ is a proper total coloring distinguishing adjacent
vertices by sums. Therefore, by Theorem~\ref{obs1} we find that $W_m \times P_n$ is $\tndi$ class II.

The remaining case is $n=3$ and $m\ge 4$. By the similar proof in Theorem~\ref{pnsm}
and the second coloring provided in Theorem~\ref{wnthm}, one can prove that $W_m \times P_3$ is $\tndi$ class I.
\end{proof}
\end{theorem}

\begin{theorem}
Let $W_m$ be a wheel of order $m \ge 3$ and $C_n$ be a cycle of order $n\ge 3$. Then
$\tndi (W_m \times C_n)= \Delta +2$.
\label{cnwm}
\begin{proof}
If $n$ be even, by Theorem~\ref{pnwm} we know
that $\tndi (W_m \times P_n) = \Delta+2$. So we have $\tndi (W_m \times C_n)= \Delta +2$
by using the same coloring $\psi$ and the edges between $(W_m \times C_n)_n$ and $(W_m \times C_n)_1$ are
colored by $m+4$.

Next, let $n=2k+1$ be odd for some $k\ge 1$. There exists a proper total
coloring distinguishing adjacent vertices by sums by
Theorem~\ref{pnwm}. For $n=1, 2, \cdots, 2k$, we use the same coloring $\psi$.
Suppose the coloring of $(W_m \times P_n)_1$ is $c$ and
$c'=c+2$. Color $(W_m \times C_n)_{2k+1}$ by the coloring $c'$.
We define a coloring $\phi$ on $W_m \times C_n$ as follows
$$ \phi ((W_m \times P_n)_i) = \left\{ \begin{array}{cl}
\psi(W_m \times C_n)& {\rm{if}} ~n=1, 2, \cdots, 2k, \\
c'((W_m \times C_n)_{2k+1}) & {\rm{if}} ~n=2k+1, \end{array}\right.$$
$$ \phi (h((W_m \times C_n)_i)) = \left\{ \begin{array}{cl}
\psi (h((W_m \times P_n)_i)) & {\rm{if}} ~i=1, 2, \cdots,2k-1,\\
1 & {\rm{if}} ~i=2k, \\
m+4 & {\rm{if}} ~i=2k+1. \end{array}\right.$$
As a result, the adjacent vertex distinguishing index by sum of $W_m \times C_n$ is $\Delta +2$.
\end{proof}
\end{theorem}

\begin{cor}
Let $S_m$ be a
star of order $m \ge 2$ and $C_n$ be a cycle of order $n\ge 3$. Then
$\tndi (S_m \times C_n)= \Delta +2$.
\label{cnsm}
\begin{proof}
$S_m \times C_n$ has two adjacent vertices of the maximum degree $\Delta=m+2$ and
the coloring $\phi$ used in Theorem~\ref{cnwm} is a proper total $m+4$ coloring
distinguishing adjacent vertices by sums of $S_m \times C_n$. By Theorem~\ref{obs1},
we find $\tndi (S_m \times C_n)= \Delta +2$.
\end{proof}
\end{cor}

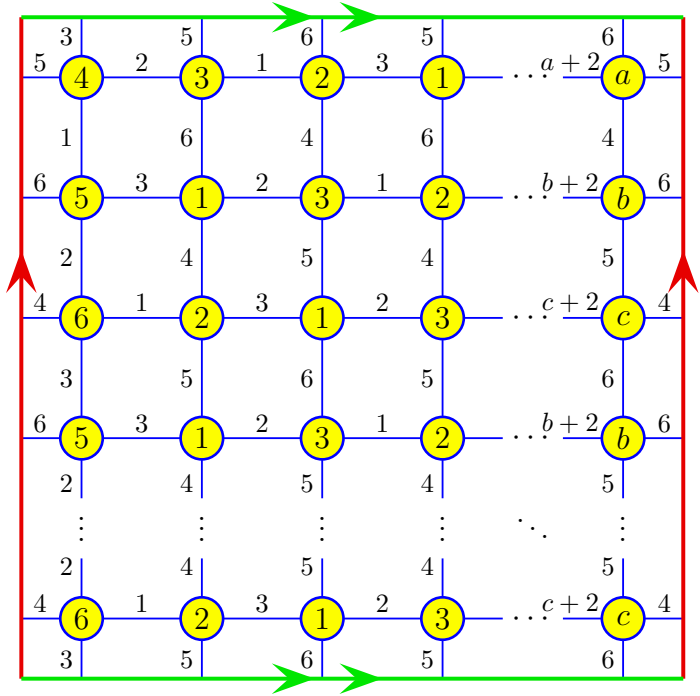
\begin{figure}
$$
\begin{pspicture}[shift=-4](-3.5,-5.9)(5.9,3.5)
\psline(-3.2,2.4)(3.2,2.4)\psline(4,2.4)(5.6,2.4)
\psline(-3.2,.8)(3.2,.8)\psline(4,.8)(5.6,.8)
\psline(-3.2,-.8)(3.2,-.8)\psline(4,-.8)(5.6,-.8)
\psline(-3.2,-2.4)(3.2,-2.4)\psline(4,-2.4)(5.6,-2.4)
\psline(-3.2,-4.8)(3.2,-4.8)\psline(4,-4.8)(5.6,-4.8)
\psline(-2.4,3.2)(-2.4,-3.2)\psline(-2.4,-4)(-2.4,-5.6)
\psline(-.8,3.2)(-.8,-3.2)\psline(-.8,-4)(-.8,-5.6)
\psline(.8,3.2)(.8,-3.2)\psline(.8,-4)(.8,-5.6)
\psline(2.4,3.2)(2.4,-3.2)\psline(2.4,-4)(2.4,-5.6)
\psline(4.8,3.2)(4.8,-3.2)\psline(4.8,-4)(4.8,-5.6)
\psline[linecolor=darkred, linewidth=1.5pt](-3.2,-5.6)(-3.2,3.2)
\psline[linecolor=darkred, arrowscale=4]{->}(-3.2,-.1)(-3.2,.1)
\psline[linecolor=darkred, linewidth=1.5pt](5.6,-5.6)(5.6,3.2)
\psline[linecolor=darkred, arrowscale=4]{->}(5.6,-.1)(5.6,.1)
\psline[linecolor=emgreen, linewidth=1.5pt](-3.2,-5.6)(5.6,-5.6)
\psline[linecolor=emgreen, arrowscale=4]{->}(.6,-5.6)(.7,-5.6)
\psline[linecolor=emgreen, arrowscale=4]{->}(1.5,-5.6)(1.6,-5.6)
\psline[linecolor=emgreen, linewidth=1.5pt](-3.2,3.2)(5.6,3.2)
\psline[linecolor=emgreen, arrowscale=4]{->}(.6,3.2)(.7,3.2)
\psline[linecolor=emgreen, arrowscale=4]{->}(1.5,3.2)(1.6,3.2)
\pscircle[fillcolor=lightgray, fillstyle=solid, linewidth=1pt](-2.4,2.4){.3}
\pscircle[fillcolor=lightgray, fillstyle=solid, linewidth=1pt](-2.4,.8){.3}
\pscircle[fillcolor=lightgray, fillstyle=solid, linewidth=1pt](-2.4,-.8){.3}
\pscircle[fillcolor=lightgray, fillstyle=solid, linewidth=1pt](-2.4,-2.4){.3}
\pscircle[fillcolor=lightgray, fillstyle=solid, linewidth=1pt](-2.4,-4.8){.3}
\pscircle[fillcolor=lightgray, fillstyle=solid, linewidth=1pt](-.8,2.4){.3}
\pscircle[fillcolor=lightgray, fillstyle=solid, linewidth=1pt](-.8,.8){.3}
\pscircle[fillcolor=lightgray, fillstyle=solid, linewidth=1pt](-.8,-.8){.3}
\pscircle[fillcolor=lightgray, fillstyle=solid, linewidth=1pt](-.8,-2.4){.3}
\pscircle[fillcolor=lightgray, fillstyle=solid, linewidth=1pt](-.8,-4.8){.3}
\pscircle[fillcolor=lightgray, fillstyle=solid, linewidth=1pt](.8,2.4){.3}
\pscircle[fillcolor=lightgray, fillstyle=solid, linewidth=1pt](.8,.8){.3}
\pscircle[fillcolor=lightgray, fillstyle=solid, linewidth=1pt](.8,-.8){.3}
\pscircle[fillcolor=lightgray, fillstyle=solid, linewidth=1pt](.8,-2.4){.3}
\pscircle[fillcolor=lightgray, fillstyle=solid, linewidth=1pt](.8,-4.8){.3}
\pscircle[fillcolor=lightgray, fillstyle=solid, linewidth=1pt](2.4,2.4){.3}
\pscircle[fillcolor=lightgray, fillstyle=solid, linewidth=1pt](2.4,.8){.3}
\pscircle[fillcolor=lightgray, fillstyle=solid, linewidth=1pt](2.4,-.8){.3}
\pscircle[fillcolor=lightgray, fillstyle=solid, linewidth=1pt](2.4,-2.4){.3}
\pscircle[fillcolor=lightgray, fillstyle=solid, linewidth=1pt](2.4,-4.8){.3}
\pscircle[fillcolor=lightgray, fillstyle=solid, linewidth=1pt](4.8,2.4){.3}
\pscircle[fillcolor=lightgray, fillstyle=solid, linewidth=1pt](4.8,.8){.3}
\pscircle[fillcolor=lightgray, fillstyle=solid, linewidth=1pt](4.8,-.8){.3}
\pscircle[fillcolor=lightgray, fillstyle=solid, linewidth=1pt](4.8,-2.4){.3}
\pscircle[fillcolor=lightgray, fillstyle=solid, linewidth=1pt](4.8,-4.8){.3}
\rput(-2.4,2.4){{$4$}} \rput(-.8,2.4){{$3$}} \rput(.8,2.4){{$2$}} \rput(2.4,2.4){{$1$}} \rput(4.8,2.4){{$a$}}
\rput(-2.4,.8){{$5$}} \rput(-.8,.8){{$1$}} \rput(.8,.8){{$3$}} \rput(2.4,.8){{$2$}} \rput(4.8,.8){{$b$}}
\rput(-2.4,-.8){{$6$}} \rput(-.8,-.8){{$2$}} \rput(.8,-.8){{$1$}} \rput(2.4,-.8){{$3$}} \rput(4.8,-.8){{$c$}}
\rput(-2.4,-2.4){{$5$}} \rput(-.8,-2.4){{$1$}} \rput(.8,-2.4){{$3$}} \rput(2.4,-2.4){{$2$}} \rput(4.8,-2.4){{$b$}}
\rput(-2.4,-4.8){{$6$}} \rput(-.8,-4.8){{$2$}} \rput(.8,-4.8){{$1$}} \rput(2.4,-4.8){{$3$}} \rput(4.8,-4.8){{$c$}}
\rput(-2.95,2.6){{\footnotesize{5}}} \rput(-1.6,2.6){{\footnotesize{2}}}
\rput(0,2.6){{\footnotesize{1}}} \rput(1.6,2.6){{\footnotesize{3}}}
\rput(4.1,2.6){{\footnotesize{$a+2$}}} \rput(5.35,2.6){{\footnotesize{5}}}
\rput(-2.95,1){{\footnotesize{6}}} \rput(-1.6,1){{\footnotesize{3}}}
\rput(0,1){{\footnotesize{2}}} \rput(1.6,1){{\footnotesize{1}}}
\rput(4.1,1){{\footnotesize{$b+2$}}} \rput(5.35,1){{\footnotesize{6}}}
\rput(-2.95,-0.6){{\footnotesize{4}}} \rput(-1.6,-0.6){{\footnotesize{1}}}
\rput(0,-0.6){{\footnotesize{3}}} \rput(1.6,-0.6){{\footnotesize{2}}}
\rput(4.1,-0.6){{\footnotesize{$c+2$}}} \rput(5.35,-0.6){{\footnotesize{4}}}
\rput(-2.95,-2.2){{\footnotesize{6}}} \rput(-1.6,-2.2){{\footnotesize{3}}}
\rput(0,-2.2){{\footnotesize{2}}} \rput(1.6,-2.2){{\footnotesize{1}}}
\rput(4.1,-2.2){{\footnotesize{$b+2$}}} \rput(5.35,-2.2){{\footnotesize{6}}}
\rput(-2.95,-4.6){{\footnotesize{4}}} \rput(-1.6,-4.6){{\footnotesize{1}}}
\rput(0,-4.6){{\footnotesize{3}}} \rput(1.6,-4.6){{\footnotesize{2}}}
\rput(4.1,-4.6){{\footnotesize{$c+2$}}} \rput(5.35,-4.6){{\footnotesize{4}}}
\rput(-2.6,2.95){{\footnotesize{3}}} \rput(-1,2.95){{\footnotesize{5}}}
\rput(0.6,2.95){{\footnotesize{6}}} \rput(2.2,2.95){{\footnotesize{5}}}
\rput(4.6,2.95){{\footnotesize{6}}}
\rput(-2.6,1.6){{\footnotesize{1}}} \rput(-1,1.6){{\footnotesize{6}}}
\rput(0.6,1.6){{\footnotesize{4}}} \rput(2.2,1.6){{\footnotesize{6}}}
\rput(4.6,1.6){{\footnotesize{4}}}
\rput(-2.6,0){{\footnotesize{2}}} \rput(-1,0){{\footnotesize{4}}}
\rput(0.6,0){{\footnotesize{5}}} \rput(2.2,0){{\footnotesize{4}}}
\rput(4.6,0){{\footnotesize{5}}}
\rput(-2.6,-1.6){{\footnotesize{3}}} \rput(-1,-1.6){{\footnotesize{5}}}
\rput(0.6,-1.6){{\footnotesize{6}}} \rput(2.2,-1.6){{\footnotesize{5}}}
\rput(4.6,-1.6){{\footnotesize{6}}}
\rput(-2.6,-3){{\footnotesize{2}}} \rput(-1,-3){{\footnotesize{4}}}
\rput(0.6,-3){{\footnotesize{5}}} \rput(2.2,-3){{\footnotesize{4}}}
\rput(4.6,-3){{\footnotesize{5}}}
\rput(-2.6,-4.1){{\footnotesize{2}}} \rput(-1,-4.1){{\footnotesize{4}}}
\rput(0.6,-4.1){{\footnotesize{5}}} \rput(2.2,-4.1){{\footnotesize{4}}}
\rput(4.6,-4.1){{\footnotesize{5}}}
\rput(-2.6,-5.35){{\footnotesize{3}}} \rput(-1,-5.35){{\footnotesize{5}}}
\rput(0.6,-5.35){{\footnotesize{6}}} \rput(2.2,-5.35){{\footnotesize{5}}}
 \rput(4.6,-5.35){{\footnotesize{6}}}
\rput(3.6,2.4){{$\cdots$}} \rput(3.6,.8){{$\cdots$}} \rput(3.6,-.8){{$\cdots$}}\rput(3.6,-2.4){{$\cdots$}} \rput(3.6,-4.8){{$\cdots$}}
\rput(-2.4,-3.5){{$\vdots$}} \rput(-.8,-3.5){{$\vdots$}} \rput(.8,-3.5){{$\vdots$}}
\rput(2.4,-3.5){{$\vdots$}} \rput(4.8,-3.5){{$\vdots$}} \rput(3.6,-3.5){{$\ddots$}}
\end{pspicture}$$
\caption{A proper total coloring distinguishing adjacent
vertices by sums of the product graph $C_n \times C_m$  where $m, n$ are odd.} \label{cncmfig9}
\end{figure}

\begin{theorem}
 Let $C_n, C_m$ be a cycle of order $n, m \ge 3$. Then
$\tndi (C_n \times C_m) = \Delta +2$.
\label{cncm}
\begin{proof}
We proceed the proof depend on the parity of $n, m$.
First suppose $n, m$ are odd. To define a coloring $\psi$ of $C_n \times C_m$,
we make a coloring $\phi$ on the first three copies
by coloring the vertices of $(C_n \times C_m)_1$ with colors $4, 3, 2, 1, 3$, $2$, $1$, $\ldots$,
and the edges of $(C_n \times C_m)_1$ with $5, 2, 1, 3, 2$, $1$, $3$, $2$, $1$, $3$, $\ldots$.
Color the vertices of $(C_n \times C_m)_2$ with colors $5, 1, 3, 2, 1, 3, 2, 1, 3, 2, \ldots$,
and the edges with $6, 3, 2, 1, 3, 2, 1, \ldots$. Color the vertices
of $(C_n \times C_m)_3$ with colors $6, 2, 1, 3, 2, 1$, $3, 2$, $1, 3, \ldots$,
and the edges with $4, 1, 3, 2, 1, 3, 2, 1, 3, 2, \ldots$. Color the edges
of $h((C_n \times C_m)_0)$ with colors $3, 5, 6, 5, 6, \ldots$, $5$, $6$, and
the edges of $h((C_n \times C_m)_1)$ with colors $1, 6, 4, 6, 4, 6, 4, \ldots, 6, 4$,
and the edges of $h((C_n \times C_m)_2)$ with colors $2, 4, 5, 4, 5, \ldots,4,5$.
Then we expend these coloring $\phi$ to $\psi$ of $C_n \times C_m$ by
$$ \psi ((C_n \times C_m)_i) = \left\{ \begin{array}{cl}
\phi((C_n \times C_m)_1)& \rm{if} ~\textit{i}~\rm{=}~\rm{1}, \\
\phi((C_n \times C_m)_2) & \rm{if} ~\textit{i}~\rm{is}~\rm{even},\\
\phi((C_n \times C_m)_3) & \rm{if} ~\textit{i}~\rm{is}~\rm{odd}~\rm{and}~\rm{bigger}~\rm{than}~\rm{3}, \end{array}\right.$$
$$ \psi (h((C_n \times C_m)_i)) = \left\{ \begin{array}{cl}
\phi(h((C_n \times C_m)_i))& {\rm{if}} ~i=0 ~{\rm{and}}~i=1, \\
\phi(h((C_n \times C_m)_2)) & {\rm{if}} ~i~{\rm{is}~\rm{even}~\rm{and}~\rm{nonzero}},\\
\phi(h((C_n \times C_m)_0)) & {\rm{if}} ~i~{\rm{is}~\rm{odd}~\rm{and}} \ge 3, \end{array}\right.$$
as illustrated in Figure~\ref{cncmfig9}
where $a\equiv 2-n$, $b\equiv 3-n$ and $c\equiv 4-n$ (mod $3$) are taken in
the standard complete residue system $\{ 1, 2, 3\}$.

\begin{figure}
$$
\begin{pspicture}[shift=-4](-3.5,-5.9)(5.9,3.5)
\psline(-3.2,2.4)(3.2,2.4)\psline(4,2.4)(5.6,2.4)
\psline(-3.2,.8)(3.2,.8)\psline(4,.8)(5.6,.8)
\psline(-3.2,-.8)(3.2,-.8)\psline(4,-.8)(5.6,-.8)
\psline(-3.2,-2.4)(3.2,-2.4)\psline(4,-2.4)(5.6,-2.4)
\psline(-3.2,-4.8)(3.2,-4.8)\psline(4,-4.8)(5.6,-4.8)
\psline(-2.4,3.2)(-2.4,-3.2)\psline(-2.4,-4)(-2.4,-5.6)
\psline(-.8,3.2)(-.8,-3.2)\psline(-.8,-4)(-.8,-5.6)
\psline(.8,3.2)(.8,-3.2)\psline(.8,-4)(.8,-5.6)
\psline(2.4,3.2)(2.4,-3.2)\psline(2.4,-4)(2.4,-5.6)
\psline(4.8,3.2)(4.8,-3.2)\psline(4.8,-4)(4.8,-5.6)
\psline[linecolor=darkred, linewidth=1.5pt](-3.2,-5.6)(-3.2,3.2)
\psline[linecolor=darkred, arrowscale=4]{->}(-3.2,-.1)(-3.2,.1)
\psline[linecolor=darkred, linewidth=1.5pt](5.6,-5.6)(5.6,3.2)
\psline[linecolor=darkred, arrowscale=4]{->}(5.6,-.1)(5.6,.1)
\psline[linecolor=emgreen, linewidth=1.5pt](-3.2,-5.6)(5.6,-5.6)
\psline[linecolor=emgreen, arrowscale=4]{->}(.6,-5.6)(.7,-5.6)
\psline[linecolor=emgreen, arrowscale=4]{->}(1.5,-5.6)(1.6,-5.6)
\psline[linecolor=emgreen, linewidth=1.5pt](-3.2,3.2)(5.6,3.2)
\psline[linecolor=emgreen, arrowscale=4]{->}(.6,3.2)(.7,3.2)
\psline[linecolor=emgreen, arrowscale=4]{->}(1.5,3.2)(1.6,3.2)
\pscircle[fillcolor=lightgray, fillstyle=solid, linewidth=1pt](-2.4,2.4){.3}
\pscircle[fillcolor=lightgray, fillstyle=solid, linewidth=1pt](-2.4,.8){.3}
\pscircle[fillcolor=lightgray, fillstyle=solid, linewidth=1pt](-2.4,-.8){.3}
\pscircle[fillcolor=lightgray, fillstyle=solid, linewidth=1pt](-2.4,-2.4){.3}
\pscircle[fillcolor=lightgray, fillstyle=solid, linewidth=1pt](-2.4,-4.8){.3}
\pscircle[fillcolor=lightgray, fillstyle=solid, linewidth=1pt](-.8,2.4){.3}
\pscircle[fillcolor=lightgray, fillstyle=solid, linewidth=1pt](-.8,.8){.3}
\pscircle[fillcolor=lightgray, fillstyle=solid, linewidth=1pt](-.8,-.8){.3}
\pscircle[fillcolor=lightgray, fillstyle=solid, linewidth=1pt](-.8,-2.4){.3}
\pscircle[fillcolor=lightgray, fillstyle=solid, linewidth=1pt](-.8,-4.8){.3}
\pscircle[fillcolor=lightgray, fillstyle=solid, linewidth=1pt](.8,2.4){.3}
\pscircle[fillcolor=lightgray, fillstyle=solid, linewidth=1pt](.8,.8){.3}
\pscircle[fillcolor=lightgray, fillstyle=solid, linewidth=1pt](.8,-.8){.3}
\pscircle[fillcolor=lightgray, fillstyle=solid, linewidth=1pt](.8,-2.4){.3}
\pscircle[fillcolor=lightgray, fillstyle=solid, linewidth=1pt](.8,-4.8){.3}
\pscircle[fillcolor=lightgray, fillstyle=solid, linewidth=1pt](2.4,2.4){.3}
\pscircle[fillcolor=lightgray, fillstyle=solid, linewidth=1pt](2.4,.8){.3}
\pscircle[fillcolor=lightgray, fillstyle=solid, linewidth=1pt](2.4,-.8){.3}
\pscircle[fillcolor=lightgray, fillstyle=solid, linewidth=1pt](2.4,-2.4){.3}
\pscircle[fillcolor=lightgray, fillstyle=solid, linewidth=1pt](2.4,-4.8){.3}
\pscircle[fillcolor=lightgray, fillstyle=solid, linewidth=1pt](4.8,2.4){.3}
\pscircle[fillcolor=lightgray, fillstyle=solid, linewidth=1pt](4.8,.8){.3}
\pscircle[fillcolor=lightgray, fillstyle=solid, linewidth=1pt](4.8,-.8){.3}
\pscircle[fillcolor=lightgray, fillstyle=solid, linewidth=1pt](4.8,-2.4){.3}
\pscircle[fillcolor=lightgray, fillstyle=solid, linewidth=1pt](4.8,-4.8){.3}
\rput(-2.4,2.4){{$1$}} \rput(-.8,2.4){{$3$}} \rput(.8,2.4){{$1$}} \rput(2.4,2.4){{$3$}} \rput(4.8,2.4){{$3$}}
\rput(-2.4,.8){{$3$}} \rput(-.8,.8){{$1$}} \rput(.8,.8){{$3$}} \rput(2.4,.8){{$1$}} \rput(4.8,.8){{$1$}}
\rput(-2.4,-.8){{$1$}} \rput(-.8,-.8){{$3$}} \rput(.8,-.8){{$1$}} \rput(2.4,-.8){{$3$}} \rput(4.8,-.8){{$3$}}
\rput(-2.4,-2.4){{$3$}} \rput(-.8,-2.4){{$1$}} \rput(.8,-2.4){{$3$}} \rput(2.4,-2.4){{$1$}} \rput(4.8,-2.4){{$1$}}
\rput(-2.4,-4.8){{$3$}} \rput(-.8,-4.8){{$1$}} \rput(.8,-4.8){{$3$}} \rput(2.4,-4.8){{$1$}} \rput(4.8,-4.8){{$1$}}
\rput(-2.95,2.6){{\footnotesize{4}}} \rput(-1.6,2.6){{\footnotesize{2}}} \rput(0,2.6){{\footnotesize{4}}} \rput(1.6,2.6){{\footnotesize{2}}} \rput(4.1,2.6){{\footnotesize{2}}} \rput(5.35,2.6){{\footnotesize{4}}}
\rput(-2.95,1){{\footnotesize{2}}} \rput(-1.6,1){{\footnotesize{4}}} \rput(0,1){{\footnotesize{2}}} \rput(1.6,1){{\footnotesize{4}}} \rput(4.1,1){{\footnotesize{4}}} \rput(5.35,1){{\footnotesize{2}}}
\rput(-2.95,-0.6){{\footnotesize{4}}} \rput(-1.6,-0.6){{\footnotesize{2}}} \rput(0,-0.6){{\footnotesize{4}}} \rput(1.6,-0.6){{\footnotesize{2}}}\rput(4.1,-0.6){{\footnotesize{2}}} \rput(5.35,-0.6){{\footnotesize{4}}}
\rput(-2.95,-2.2){{\footnotesize{2}}} \rput(-1.6,-2.2){{\footnotesize{4}}} \rput(0,-2.2){{\footnotesize{2}}} \rput(1.6,-2.2){{\footnotesize{4}}} \rput(4.1,-2.2){{\footnotesize{4}}} \rput(5.35,-2.2){{\footnotesize{2}}}
\rput(-2.95,-4.6){{\footnotesize{2}}} \rput(-1.6,-4.6){{\footnotesize{4}}} \rput(0,-4.6){{\footnotesize{2}}} \rput(1.6,-4.6){{\footnotesize{4}}} \rput(4.1,-4.6){{\footnotesize{4}}} \rput(5.35,-4.6){{\footnotesize{2}}}
\rput(-2.6,2.95){{\footnotesize{6}}} \rput(-1,2.95){{\footnotesize{6}}} \rput(0.6,2.95){{\footnotesize{6}}} \rput(2.2,2.95){{\footnotesize{6}}} \rput(4.6,2.95){{\footnotesize{6}}}
\rput(-2.6,1.6){{\footnotesize{5}}} \rput(-1,1.6){{\footnotesize{5}}} \rput(0.6,1.6){{\footnotesize{5}}} \rput(2.2,1.6){{\footnotesize{5}}} \rput(4.6,1.6){{\footnotesize{5}}}
\rput(-2.6,0){{\footnotesize{6}}} \rput(-1,0){{\footnotesize{6}}} \rput(0.6,0){{\footnotesize{6}}} \rput(2.2,0){{\footnotesize{6}}} \rput(4.6,0){{\footnotesize{6}}}
\rput(-2.6,-1.6){{\footnotesize{5}}} \rput(-1,-1.6){{\footnotesize{5}}} \rput(0.6,-1.6){{\footnotesize{5}}} \rput(2.2,-1.6){{\footnotesize{5}}} \rput(4.6,-1.6){{\footnotesize{5}}}
\rput(-2.6,-3){{\footnotesize{6}}} \rput(-1,-3){{\footnotesize{6}}} \rput(0.6,-3){{\footnotesize{6}}} \rput(2.2,-3){{\footnotesize{6}}}  \rput(4.6,-3){{\footnotesize{6}}}
\rput(-2.6,-4.1){{\footnotesize{5}}} \rput(-1,-4.1){{\footnotesize{5}}} \rput(0.6,-4.1){{\footnotesize{5}}} \rput(2.2,-4.1){{\footnotesize{5}}} \rput(4.6,-4.1){{\footnotesize{5}}}
\rput(-2.6,-5.35){{\footnotesize{6}}} \rput(-1,-5.35){{\footnotesize{6}}} \rput(0.6,-5.35){{\footnotesize{6}}} \rput(2.2,-5.35){{\footnotesize{6}}} \rput(4.6,-5.35){{\footnotesize{6}}}
\rput(3.6,2.4){{$\cdots$}} \rput(3.6,.8){{$\cdots$}} \rput(3.6,-.8){{$\cdots$}}\rput(3.6,-2.4){{$\cdots$}} \rput(3.6,-4.8){{$\cdots$}}
\rput(-2.4,-3.5){{$\vdots$}} \rput(-.8,-3.5){{$\vdots$}} \rput(.8,-3.5){{$\vdots$}}
\rput(2.4,-3.5){{$\vdots$}} \rput(4.8,-3.5){{$\vdots$}} \rput(3.6,-3.5){{$\ddots$}}
\end{pspicture}$$
\caption{A proper total coloring distinguishing adjacent vertices by
sums of the product graph $C_n \times C_m$  where $m, n$ are even.} \label{cncmfig10}
\end{figure}
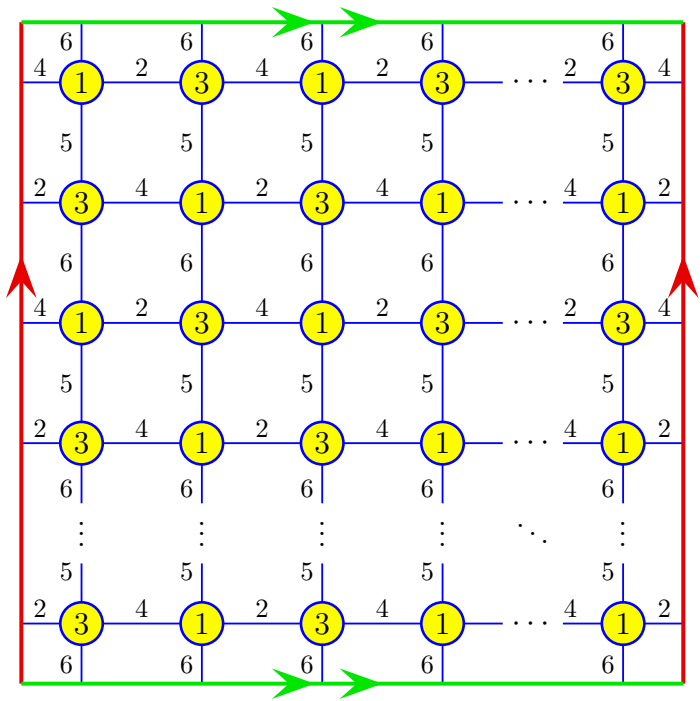

For the rest of cases, we first make a coloring $\phi$ on some copies, then expand it to
the coloring $\psi$ of $C_n \times C_m$ as we did for the case $n, m$ are odd.

Let $n, m$ be even. We color the vertices of $(C_n \times C_m)_1$
with colors $1, 3, 1, 3, \ldots, 1, 3$, and the edges with $4, 2, 4, 2, \ldots, 4, 2$.
Color the vertices of $(C_n \times C_m)_2$ with colors $3, 1, 3, 1, \ldots, 3, 1$,
and the edges with $2, 4, 2, 4, \ldots, 2, 4$.
$$ \psi ((C_n \times C_m)_i) = \left\{ \begin{array}{cl}
\phi((C_n \times C_m)_1)& \rm{if} ~\textit{i}~\rm{is}~\rm{odd}, \\
\phi((C_n \times C_m)_2) & \rm{if} ~\textit{i}~\rm{is}~\rm{even}, \end{array}\right.$$
$$ \psi (h((C_n \times C_m)_i)) = \left\{ \begin{array}{cl}
5& \rm{if} ~\textit{i}~\rm{is}~\rm{odd}, \\
6 & \rm{if} ~\textit{i}~\rm{is}~\rm{even}, \end{array}\right.$$
as illustrated in Figure~\ref{cncmfig10}.

\begin{figure}
$$
\begin{pspicture}[shift=-4](-3.5,-5.9)(5.9,3.5)
\psline(-3.2,2.4)(3.2,2.4)\psline(4,2.4)(5.6,2.4)
\psline(-3.2,.8)(3.2,.8)\psline(4,.8)(5.6,.8)
\psline(-3.2,-.8)(3.2,-.8)\psline(4,-.8)(5.6,-.8)
\psline(-3.2,-2.4)(3.2,-2.4)\psline(4,-2.4)(5.6,-2.4)
\psline(-3.2,-4.8)(3.2,-4.8)\psline(4,-4.8)(5.6,-4.8)
\psline(-2.4,3.2)(-2.4,-3.2)\psline(-2.4,-4)(-2.4,-5.6)
\psline(-.8,3.2)(-.8,-3.2)\psline(-.8,-4)(-.8,-5.6)
\psline(.8,3.2)(.8,-3.2)\psline(.8,-4)(.8,-5.6)
\psline(2.4,3.2)(2.4,-3.2)\psline(2.4,-4)(2.4,-5.6)
\psline(4.8,3.2)(4.8,-3.2)\psline(4.8,-4)(4.8,-5.6)
\psline[linecolor=darkred, linewidth=1.5pt](-3.2,-5.6)(-3.2,3.2)
\psline[linecolor=darkred, arrowscale=4]{->}(-3.2,-.1)(-3.2,.1)
\psline[linecolor=darkred, linewidth=1.5pt](5.6,-5.6)(5.6,3.2)
\psline[linecolor=darkred, arrowscale=4]{->}(5.6,-.1)(5.6,.1)
\psline[linecolor=emgreen, linewidth=1.5pt](-3.2,-5.6)(5.6,-5.6)
\psline[linecolor=emgreen, arrowscale=4]{->}(.6,-5.6)(.7,-5.6)
\psline[linecolor=emgreen, arrowscale=4]{->}(1.5,-5.6)(1.6,-5.6)
\psline[linecolor=emgreen, linewidth=1.5pt](-3.2,3.2)(5.6,3.2)
\psline[linecolor=emgreen, arrowscale=4]{->}(.6,3.2)(.7,3.2)
\psline[linecolor=emgreen, arrowscale=4]{->}(1.5,3.2)(1.6,3.2)
\pscircle[fillcolor=lightgray, fillstyle=solid, linewidth=1pt](-2.4,2.4){.3}
\pscircle[fillcolor=lightgray, fillstyle=solid, linewidth=1pt](-2.4,.8){.3}
\pscircle[fillcolor=lightgray, fillstyle=solid, linewidth=1pt](-2.4,-.8){.3}
\pscircle[fillcolor=lightgray, fillstyle=solid, linewidth=1pt](-2.4,-2.4){.3}
\pscircle[fillcolor=lightgray, fillstyle=solid, linewidth=1pt](-2.4,-4.8){.3}
\pscircle[fillcolor=lightgray, fillstyle=solid, linewidth=1pt](-.8,2.4){.3}
\pscircle[fillcolor=lightgray, fillstyle=solid, linewidth=1pt](-.8,.8){.3}
\pscircle[fillcolor=lightgray, fillstyle=solid, linewidth=1pt](-.8,-.8){.3}
\pscircle[fillcolor=lightgray, fillstyle=solid, linewidth=1pt](-.8,-2.4){.3}
\pscircle[fillcolor=lightgray, fillstyle=solid, linewidth=1pt](-.8,-4.8){.3}
\pscircle[fillcolor=lightgray, fillstyle=solid, linewidth=1pt](.8,2.4){.3}
\pscircle[fillcolor=lightgray, fillstyle=solid, linewidth=1pt](.8,.8){.3}
\pscircle[fillcolor=lightgray, fillstyle=solid, linewidth=1pt](.8,-.8){.3}
\pscircle[fillcolor=lightgray, fillstyle=solid, linewidth=1pt](.8,-2.4){.3}
\pscircle[fillcolor=lightgray, fillstyle=solid, linewidth=1pt](.8,-4.8){.3}
\pscircle[fillcolor=lightgray, fillstyle=solid, linewidth=1pt](2.4,2.4){.3}
\pscircle[fillcolor=lightgray, fillstyle=solid, linewidth=1pt](2.4,.8){.3}
\pscircle[fillcolor=lightgray, fillstyle=solid, linewidth=1pt](2.4,-.8){.3}
\pscircle[fillcolor=lightgray, fillstyle=solid, linewidth=1pt](2.4,-2.4){.3}
\pscircle[fillcolor=lightgray, fillstyle=solid, linewidth=1pt](2.4,-4.8){.3}
\pscircle[fillcolor=lightgray, fillstyle=solid, linewidth=1pt](4.8,2.4){.3}
\pscircle[fillcolor=lightgray, fillstyle=solid, linewidth=1pt](4.8,.8){.3}
\pscircle[fillcolor=lightgray, fillstyle=solid, linewidth=1pt](4.8,-.8){.3}
\pscircle[fillcolor=lightgray, fillstyle=solid, linewidth=1pt](4.8,-2.4){.3}
\pscircle[fillcolor=lightgray, fillstyle=solid, linewidth=1pt](4.8,-4.8){.3}
\rput(-2.4,2.4){{$1$}} \rput(-.8,2.4){{$3$}} \rput(.8,2.4){{$1$}} \rput(2.4,2.4){{$3$}} \rput(4.8,2.4){{$3$}}
\rput(-2.4,.8){{$3$}} \rput(-.8,.8){{$1$}} \rput(.8,.8){{$3$}} \rput(2.4,.8){{$1$}} \rput(4.8,.8){{$1$}}
\rput(-2.4,-.8){{$1$}} \rput(-.8,-.8){{$3$}} \rput(.8,-.8){{$1$}} \rput(2.4,-.8){{$3$}} \rput(4.8,-.8){{$3$}}
\rput(-2.4,-2.4){{$3$}} \rput(-.8,-2.4){{$1$}} \rput(.8,-2.4){{$3$}} \rput(2.4,-2.4){{$1$}} \rput(4.8,-2.4){{$1$}}
\rput(-2.4,-4.8){{$3$}} \rput(-.8,-4.8){{$1$}} \rput(.8,-4.8){{$3$}} \rput(2.4,-4.8){{$1$}} \rput(4.8,-4.8){{$1$}}
\rput(-2.95,2.6){{\footnotesize{4}}} \rput(-1.6,2.6){{\footnotesize{2}}} \rput(0,2.6){{\footnotesize{4}}} \rput(1.6,2.6){{\footnotesize{2}}} \rput(4.1,2.6){{\footnotesize{2}}} \rput(5.35,2.6){{\footnotesize{4}}}
\rput(-2.95,1){{\footnotesize{2}}} \rput(-1.6,1){{\footnotesize{4}}} \rput(0,1){{\footnotesize{2}}} \rput(1.6,1){{\footnotesize{4}}} \rput(4.1,1){{\footnotesize{4}}} \rput(5.35,1){{\footnotesize{2}}}
\rput(-2.95,-0.6){{\footnotesize{4}}} \rput(-1.6,-0.6){{\footnotesize{2}}} \rput(0,-0.6){{\footnotesize{4}}} \rput(1.6,-0.6){{\footnotesize{2}}}\rput(4.1,-0.6){{\footnotesize{2}}} \rput(5.35,-0.6){{\footnotesize{4}}}
\rput(-2.95,-2.2){{\footnotesize{2}}} \rput(-1.6,-2.2){{\footnotesize{4}}} \rput(0,-2.2){{\footnotesize{2}}} \rput(1.6,-2.2){{\footnotesize{4}}} \rput(4.1,-2.2){{\footnotesize{4}}} \rput(5.35,-2.2){{\footnotesize{2}}}
\rput(-2.95,-4.6){{\footnotesize{2}}} \rput(-1.6,-4.6){{\footnotesize{4}}} \rput(0,-4.6){{\footnotesize{2}}} \rput(1.6,-4.6){{\footnotesize{4}}} \rput(4.1,-4.6){{\footnotesize{4}}} \rput(5.35,-4.6){{\footnotesize{2}}}
\rput(-2.6,2.95){{\footnotesize{6}}} \rput(-1,2.95){{\footnotesize{6}}} \rput(0.6,2.95){{\footnotesize{6}}} \rput(2.2,2.95){{\footnotesize{6}}} \rput(4.6,2.95){{\footnotesize{6}}}
\rput(-2.6,1.6){{\footnotesize{5}}} \rput(-1,1.6){{\footnotesize{5}}} \rput(0.6,1.6){{\footnotesize{5}}} \rput(2.2,1.6){{\footnotesize{5}}} \rput(4.6,1.6){{\footnotesize{5}}}
\rput(-2.6,0){{\footnotesize{6}}} \rput(-1,0){{\footnotesize{6}}} \rput(0.6,0){{\footnotesize{6}}} \rput(2.2,0){{\footnotesize{6}}} \rput(4.6,0){{\footnotesize{6}}}
\rput(-2.6,-1.6){{\footnotesize{5}}} \rput(-1,-1.6){{\footnotesize{5}}} \rput(0.6,-1.6){{\footnotesize{5}}} \rput(2.2,-1.6){{\footnotesize{5}}} \rput(4.6,-1.6){{\footnotesize{5}}}
\rput(-2.6,-3){{\footnotesize{6}}} \rput(-1,-3){{\footnotesize{6}}} \rput(0.6,-3){{\footnotesize{6}}} \rput(2.2,-3){{\footnotesize{6}}}  \rput(4.6,-3){{\footnotesize{6}}}
\rput(-2.6,-4.1){{\footnotesize{5}}} \rput(-1,-4.1){{\footnotesize{5}}} \rput(0.6,-4.1){{\footnotesize{5}}} \rput(2.2,-4.1){{\footnotesize{5}}} \rput(4.6,-4.1){{\footnotesize{5}}}
\rput(-2.6,-5.35){{\footnotesize{6}}} \rput(-1,-5.35){{\footnotesize{6}}} \rput(0.6,-5.35){{\footnotesize{6}}} \rput(2.2,-5.35){{\footnotesize{6}}} \rput(4.6,-5.35){{\footnotesize{6}}}
\rput(3.6,2.4){{$\cdots$}} \rput(3.6,.8){{$\cdots$}} \rput(3.6,-.8){{$\cdots$}}\rput(3.6,-2.4){{$\cdots$}} \rput(3.6,-4.8){{$\cdots$}}
\rput(-2.4,-3.5){{$\vdots$}} \rput(-.8,-3.5){{$\vdots$}} \rput(.8,-3.5){{$\vdots$}}
\rput(2.4,-3.5){{$\vdots$}} \rput(4.8,-3.5){{$\vdots$}} \rput(3.6,-3.5){{$\ddots$}}
\end{pspicture}$$
\caption{A proper total coloring distinguishing adjacent
vertices by sums of the product graph $C_n \times C_m$
where one of $n, m$ is even.} \label{cncmfig11}
\end{figure}
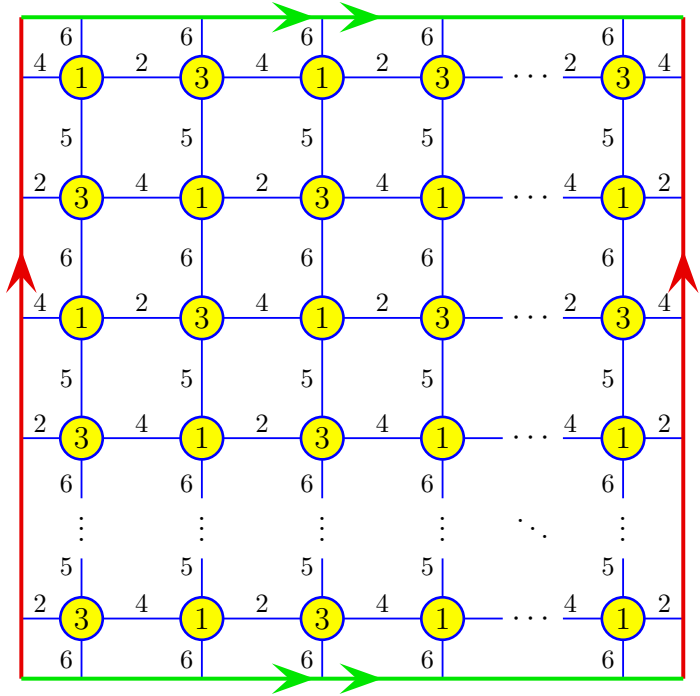

If one of $n, m$ is even and the other is odd.
By exchanging $n, m$, we my assume that $m$ is even.
If $n=3$, we use the coloring of $C_3 \times P_m$ as depicted in Figure~\ref{c3pnfig}
since $m$ is even, we can naturally extend this coloring to a coloring of $C_3 \times C_m$
and it is a proper total coloring distinguishing adjacent vertices by sums. Now we assume $n\ge 5$.
We color the vertices of $(C_n \times C_m)_1$
 with colors $4, 1, 3, 1, 3, \ldots, 1, 3$, and the edges with
 $1, 3, 2, 4, 2, 4, \ldots, 4, 2$.
 Color the vertices of $(C_n \times C_m)_2$
 with colors $2, 3, 1, 3, 1, \ldots, 3, 1$,
 and the edges with $3, 1, 4, 2, 4, \ldots, 2, 4$.
$$ \psi ((C_n \times C_m)_i) = \left\{ \begin{array}{cl}
\phi((C_n \times C_m)_1)& \rm{if} ~\textit{i}~\rm{is}~\rm{odd}, \\
\phi((C_n \times C_m)_2) & \rm{if} ~\textit{i}~\rm{is}~\rm{even}, \end{array}\right.$$
$$ \psi (h((C_n \times C_m)_i)) = \left\{ \begin{array}{cl}
5& \rm{if} ~\textit{i}~\rm{is}~\rm{odd}, \\
6 & \rm{if} ~\textit{i}~\rm{is}~\rm{even},\end{array}\right.$$
as illustrated in Figure~\ref{cncmfig11}.
These graphs always have two adjacent vertices $x, y$
such that $deg(x)=deg(y)=\Delta$. By Theorem~\ref{obs1},
we get $ \tndi (C_n \times C_m) = \Delta +2$.
\end{proof}
\end{theorem}

\section{Conclusion and discussion} \label{endcom}

As we have seen the results in Section~\ref{prod},
the most of cases, the product graphs are $\tndi$ Class II.
Although we are only able to find the
adjacent vertex distinguishing index of a few product graphs for which one of
component of the product is a path $P_n$, we expect
that the adjacent vertex distinguishing index of $\Gamma \times P_n$
can be found where $\Gamma$ is a regular graph.

It has been more than three years since a proper total colorings distinguishing
adjacent vertices by sums was first invented, however, it is still unknown that the
adjacent vertex distinguishing index by sum and adjacent vertex distinguishing index
are different or not. If these two are the same, then this will benefit the computer
program to find adjacent vertex distinguishing index because comparing the total sum is
immensely faster than comparing two sets for all adjacent vertices.
However, we know the following inequality between these two indices as follows.

Since $f(v)= \sum_{\alpha \in C(v)} \alpha$ and the implication, if $f(u) \neq f(v)$, then
$C(u) \neq C(v)$, one can easily see that if a graph $\Gamma$ has a
proper total $k$-colorings distinguishing adjacent vertices by sums,
then $\Gamma$ has an AVD total $k$-coloring. Therefore,
$$\chi_{at}(\Gamma) \le \tndi (\Gamma).$$

At last, we state two conjectures regrading $\tndi$ classes.

\begin{conj} \label{conj1}
A graph $\Gamma$ is $\tndi$ class III if and only if
$\Gamma =K_{2n+1}$ for some $n \ge 1$.
\end{conj}

Is is fairly easy to see that if a graph $\Gamma$ is regular with valency $k$ and it is $\tndi$ class II,
then at each vertex the set of colors of the vertex and its incident edges have a cardinality $k+1$ out of $k+2$
colors as demonstrated in the complete graph $K_{2n}$. Thus if we assign the missing color to the vertex, we get a proper $k+2$ coloring.
We expect the converse also can be proven.

\begin{conj} \label{conj2}
A $k$ regular graph $\Gamma$ is $\tndi$ class II if and only if
$\Gamma$ has a proper vertex $k+2$ coloring.
\end{conj}

\section*{Acknowledgments}
The \TeX\, macro package
PSTricks~\cite{PSTricks} was essential for typesetting the equations
and figures. This work was supported by the Korea Foundation for the Advancement of Science \& Creativity(KOFAC), and funded by the Korean Government(MOE).

\end{document}